\documentclass{siamonline190516}
\usepackage{latex-macros}
\usepackage[capitalize]{cleveref}
\usepackage{algpseudocode}
\usepackage{appendix}
\usepackage{cancel}
\usepackage{enumitem}
\usepackage{booktabs}

\definecolor{lblue}{HTML}{908cc0}
\definecolor{mblue}{HTML}{519cc8}
\definecolor{hblue}{HTML}{1d5996}
\definecolor{lred}{HTML}{cb5501}
\definecolor{mred}{HTML}{f1885b}
\definecolor{hred}{HTML}{b3001e}
\definecolor{ttred}{HTML}{ca3542}

\crefname{appsec}{Appendix}{Appendices}

\newtheorem{assumption}{Assumption}
\newtheorem{remark}{Remark}

\newcommand{\tsup}[1]{\textsuperscript{#1}}

\title{Incrementally Updated Spectral Embeddings}
\author{
Vasileios Charisopoulos\thanks{Department of Operations
Research \& Information Engineering, Cornell University, 14850 Ithaca, NY.
Email: \texttt{vc333@cornell.edu}}
\and
Austin R.\ Benson\thanks{Department of Computer Science, Cornell University, 14850
Ithaca, NY. Email: \texttt{arb@cs.cornell.edu}}
\and
Anil Damle\thanks{Department of Computer Science,
Cornell University, 14850 Ithaca, NY. Email: \texttt{damle@cornell.edu}}}

\newcommand{\diag}{\mathrm{diag}}

\newcommand{\mxhat}[2]{{\hat{{#1}}_{({#2})}}}

\begin{document}

\maketitle

\begin{abstract}
    Several fundamental tasks in data science rely on computing
    an extremal eigenspace of size $r \ll n$, where $n$ is the underlying
    problem dimension. For example, spectral clustering and PCA both
    require the computation of the leading $r$-dimensional subspace.
    Often, this process is repeated
    over time due to the possible temporal nature of
    the data; \emph{e.g.}, graphs representing relations in a social network may
    change over time, and feature vectors may be added, removed
    or updated in a dataset. Therefore, it is important to efficiently carry out
    the computations involved to keep up with frequent changes
    in the underlying data and also to dynamically determine a reasonable
    size for the subspace of interest.
    We present a complete computational pipeline for efficiently updating
    spectral embeddings in a variety of contexts.
    Our basic approach is to ``seed'' iterative methods for
    eigenproblems with the most recent subspace estimate to significantly
    reduce the computations involved, in contrast with a na\"ive approach which
    recomputes the subspace of interest from scratch at every step.
    In this setting, we provide various bounds on the number of iterations
    common eigensolvers need to perform in order to update the extremal
    eigenspace to a sufficient tolerance. We also incorporate a criterion for
    determining the size of the subspace based on successive eigenvalue
    ratios.
    We demonstrate the merits of our approach on the tasks of spectral
    clustering of temporally evolving graphs and PCA of an incrementally
    updated data matrix.
\end{abstract}



\section{Introduction}
In the big data era, scientists and engineers need to operate on
massive datasets on a daily basis, fueling essential algorithms for
commercial or scientific applications. These datasets can contain millions or
billions of data points, making the task of extracting meaningful
information especially challenging~\cite{LesRajUll14}; moreover, it is often
the case that the data are also high-dimensional, which can significantly affect
the time and storage required to work with such datasets. These computational
challenges strongly motivate the need to work with ``summarized''
versions of data that facilitate fast and memory-efficient computations.
Several popular techniques in data science have been devised to address this
problem, such as sketching~\cite{Woodruff14}, dimensionality
reduction~\cite{LDA98,Jolliffe02,LLE00}, and limited-memory/stochastic
algorithms in optimization (see, \emph{e.g.,}~\cite{BotCurNoc18} and references therein).

A popular approach to remedy some of the aforementioned difficulties
revolves around using \textbf{spectral} information associated to the problem.
For example, many high-dimensional problems involve data
that can be modeled as \textit{graphs} encoding interactions between
entities in social networks~\cite{PanOpsCar09}, gene interaction
data~\cite{OVER+08} and product recommendation networks~\cite{LesAdamHub07}. A
task of fundamental importance is clustering the nodes of the network
into communities or clusters of similar
nodes~\cite{Porter-2009-communities,Schaeffer-2007-clustering}.
\textit{Spectral clustering} is a popular approach for this problem that computes
the $r$ leading eigenvectors of the (symmetrically normalized)
adjacency matrix, using this subspace as a low-dimensional
representation of nodes in the graph, and then feeding this representation to
a point cloud clustering algorithm such as \texttt{k-means}~\cite{NgJorWei02}.
Determining $r$ on the fly is often done by comparing the
ratios of successive eigenvalues~\cite{VonLux07}.

In a similar vein, linear dimensionality reduction is often tackled via
Principal Component Analysis (PCA)~\cite{Jolliffe02}, which projects the
data matrix $X$ to a low-dimensional coordinate system which is
spanned by the leading eigenvectors of the covariance matrix, $X^\top X$. The
$i$\tsup{th} principal component is precisely the projection of the data matrix $X$
to the $i$\tsup{th} eigenvector. In nonlinear dimensionality reduction, a common
assumption is that the data lie on a low-dimensional manifold embedded in
high-dimensional space. The method of Laplacian eigenmaps~\cite{BelkinNiyogi03},
which operates under this assumption, constructs a weighted adjacency matrix
and solves a small sequence of eigenproblems to compute a lower-dimensional
embedding of the data points.

\subsection{Dealing with data that gets updated} \label{sec:intro-temporal}
Going one step further, many datasets get incrementally updated over time.
For example, graph data can change in many ways:
new links form between entities,
old connections cease to exist, and
new nodes are introduced to a network~\cite{Viswanath-2009-evolution,Kumar-2010-structure}.
In such cases, old clustering assignments may no longer
accurately reflect the community structure of the network, and therefore need
to be updated~\cite{NXCG+07}. In other cases, the dimensionality of data may
increase as more informative features as well as new data become
available~\cite{BalDasFreu13,LJBL+11}. Such updates reduce the fidelity of the previously
computed low-dimensional embeddings. In both
cases, we are faced with a pressing question: \textit{has the quality of the
spectral embedding degraded, and if so, can we update it efficiently?}

The main theoretical tools for answering the first part of this question come
from matrix perturbation theory~\cite{StewSun90}, which provides worst-case
bounds for the distance between eigenspaces. However, algorithms proposed for
updating spectral embeddings may not necessarily utilize those results. For
example, if the original (unperturbed) data matrix \emph{and} its updates are low rank,
then one can update the thin SVD in practically linear time~\cite{Brand06};
however, this is rarely the case in practice.
Other ``direct'' approaches under low-rank updates are often inapplicable since
they require the full set of singular values and vectors to be available for the
unperturbed matrix~\cite{BunchNielsen78,GuEisenstat96}.
When the data matrix is sparse, a natural candidate may be
updating the sparse Cholesky or $L D L^\top$ decompositions~\cite{Cholmod08};
however, these rely on heuristics and may not necessarily preserve sparsity
across multiple updates and it is unclear how to adapt these methods
to maintain a low-dimensional subspace.

Several of the aforementioned approaches share a common prohibitive requirement: the
knowledge of the full set of eigenvalues and eigenvectors at the beginning of an
update step. If one is willing to forego guarantees on the accuracy of the
updated embedding, it is at least empirically possible
to exploit the first-order expansion of eigenvalues, where explicit
knowledge of the full set is no longer required~\cite{NXCG+07}.
On the downside, this entails significant uncertainty about the quality of the maintained embedding.
A more promising and relevant heuristic is that of \textit{warm-starting}
iterative eigenvalue methods. Many eigensolvers are typically initialized using a random
guess, but when the target is a sequence of similar or related eigenproblems it is
natural to expect that the previous solution might be a good initial estimate.
This heuristic has been applied towards solving sequences of linear
systems~\cite{KilmerDesturier06,PDMJ+06}, sometimes referred to as
\textit{recycled Krylov method}, and problems arising from successive
linearizations of nonlinear eigenvalue problems~\cite{SalGirSaadMor15}.

Somewhat surprisingly, despite the empirical success of warm starting, it
is often overlooked in favor of other heuristics~\cite{NXCG+07} or matrix
sketching to reduce per-iteration cost~\cite{Gu15,HalMarTropp11}.

\subsection{Our contribution}
Given the preceding discussion, an ideal method for
incrementally computing spectral embeddings should:
\begin{enumerate}
	\item incur minimal additional work per modification, ideally independent
    of the underlying problem's dimension; \label{point:1}
	\item require minimal additional storage or extensive pre-computation (this excludes
	         classical approaches that assume knowledge of the full set of eigenvalues~\cite{BunchNielsen78}); and
	\item take advantage of any features of the underlying matrix, such as
		sparsity or structure (e.g. Toeplitz matrices).
\end{enumerate}
With regards to~\cref{point:1}, it is important to be able to bound the amount
of work required to update a spectral embedding \textit{on the fly}, before
performing the actual update.

In this paper, we show how warm-started iterative methods can naturally satisfy all of
the above criteria under a minimal assumption on eigenvalue decay formalized
in \Cref{sec:time-subspace-estimation}. We show that when the perturbations are
``small'', using the previous subspace estimate as a ``seed'' for warm-starting entails
essentially $\cO(1)$ additional iterations per update even when the previous
estimates have been computed inexactly (orders of magnitude above floating point
error). We provide an adaptive algorithm for updating the spectral embedding
over time as well \textit{computable bounds} on the number of iterations
required after each modification to the underlying matrix, for two standard
iterative methods (eigensolvers).
Because the main bottleneck operation of iterative methods is usually matrix-vector
multiplication~\cite{Demmel97,GVL13,Saad11}, the proposed algorithm naturally
accelerates whenever the underlying matrix is structured or sparse. It also
enables tracking a proxy for the appropriate subspace dimension, which can
aid in determining the correct size of the spectral embedding on the fly.

In \Cref{sec:applications}, we present a few concrete applications which
fit into the incremental framework outlined above. In these settings, we are
able to utilize properties of the perturbations to derive further a priori
bounds on their effects on subspace distance. For example, when the
perturbations are outer products of Gaussian random vectors, the subspace distance
is affected by a factor of at most $\tilde{\cO}\left(\sqrt{\frac{r}{n}}\right)$
with high probability. We also present experiments on real and synthetic data
which validate the effectiveness of our approach.

Finally, we stress that we do not intend to provide a replacement to or compete
with established eigensolvers, but rather describe a concrete pipeline for
incremental spectral embeddings with well-understood worst-case guarantees.
Indeed, our method only benefits from advances in software aimed at solving
eigenproblems.

\subsection{Notation}
Throughout the paper, we endow $\Rbb^{n}$ with the standard inner product
$\ip{x, y} := x^\top y$ and the induced norm $\norm{x}_2 = \sqrt{\ip{x, x}}$.
Given a matrix $A \in \Rbb^{n_1 \times n_2}$ , we
write $\norm{A}_2$ for its spectral norm $\sup_{\norm{x}_2 = 1} \norm{Ax}_2$
and $\norm{A}_F := \sqrt{\mathrm{tr}(A^\top A)}$ for its Frobenius norm.
We let $\Obb_{n, k} := \set{Q \in \Rbb^{n \times k} \mmid Q^\top Q = I}$
denote the set of real $n \times k$ matrices with orthonormal columns and drop the second
subscript when $k = n$.

We write $\lambda_i(A)$ for the $i$\tsup{th} eigenvalue of a matrix $A$, assuming a
descending order such that $\lambda_1(A) \geq \dots \geq \lambda_n(A)$.
We follow the same notation for singular values,
$\sigma_i(A)$. Given a symmetric matrix $A \in \Rbb^{n \times n}$, we will use
its \textbf{spectral decomposition}:
\begin{equation}
	A = \begin{bmatrix}
		V & V_{\perp}
	\end{bmatrix} \begin{bmatrix}
		\Lambda & 0 \\ 0 & \Lambda_{\perp}
	\end{bmatrix}
	\begin{bmatrix}
		V^\top \\
		V_{\perp}^\top
	\end{bmatrix},
    \label{eq:schur-decomp}
\end{equation}
where $\bmx{V & V_{\perp}}$ are the eigenvectors of $A$ and $\bmx{\Lambda & 0 \\
0 & \Lambda_{\perp}}$ is a diagonal matrix containing the eigenvalues of $A$.
By convention, $V$ should be understood to correspond to the subspace of interest.
Similarly, when $A \in \Rbb^{m \times n}$ with $m \geq n$ is an arbitrary rectangular matrix,
we will use its \textbf{singular value decomposition}:
\[
	A = U \Sigma V^\top
\]
where $U \in \Obb_{m, n}, V \in \Obb_{n}$ are the left and right singular vectors of $A$
and $\Sigma$ is a diagonal matrix containing the singular values of $A$. We
will write $\kappa(A) := \frac{\sigma_1(A)}{\sigma_n(A)}$ for $A$'s condition
number. Throughout the main text, we refer to standard linear algebraic
routines, which
are built into numerical software packages such as \texttt{LAPACK}, using the
following notation:
\begin{itemize}
\item $Q, R = \texttt{qr}(A)$ for the (reduced) QR factorization of $A \in
\Rbb^{m \times n}$ with $m\geq n$, which decomposes $A$ into an orthogonal matrix
$Q \in \Obb_{m, n}$ and an upper triangular matrix $R \in \Rbb^{n \times n}$.
\item $V, \Lambda = \texttt{spectral}(A)$ for the spectral decomposition of a
symmetric matrix $A \in \Rbb^{n \times n}$, with $V$ containing the
eigenvectors and $\Lambda := \mathrm{diag}(\lambda_1, \dots, \lambda_n)$ containing the
eigenvalues of $A$, respectively. We will assume that $\lambda_1 \geq
\dots \geq \lambda_n$ in $\Lambda$.
\end{itemize}
Moreover, we use $\diag(x_1, \dots, x_n)$ to refer to the diagonal matrix
constructed from the vector $x = \begin{pmatrix} x_1 & \dots & x_n \end{pmatrix}^\top$, and when
$D$ is a matrix denote $\diag(D) := \begin{pmatrix} D_{11} & D_{22} & \dots &
D_{nn} \end{pmatrix}^\top$.
We also use notation from Golub and Van Loan~\cite{GVL13} for matrix or vector
slicing: $A_{:, 1:r}$ denotes the submatrix of $A$ formed by taking the first $r$ columns,
while $x_{1:r}$ denotes the vector formed by the leading $r$ elements of $x$.

\section{Iterative Methods for Eigenvalue Problems} \label{sec:background}

%

In this section we present two popular iterative methods for computing the
eigenvalues of symmetric / Hermitian matrices that form the building blocks of
our algorithms. Additionally, we review their convergence guarantees.
Throughout, we assume that we are interested in the subspace corresponding to
the \textbf{algebraically largest} eigenvalues (which are also the largest in
magnitude, up to a shift).
Recall that the \textit{Ritz values} are the eigenvalues of the matrix
${V_0}^\top A V_0$, where $V_0$ is an approximation to an invariant subspace
of $A$. Below, we are interested in subspace distances measured in the spectral
norm,
\(
	\mathrm{dist}(V, \tilde{V}) := \norm{VV^\top - \tilde{V} \tilde{V}^\top}_2,
\)
which has been the focus of the convergence analysis of iterative methods for
eigenproblems, and pairs well with established matrix perturbation theory for
unitarily invariant norms.
\begin{remark}
While our work predominantly relies on standard notions of subspace distance, many applications may benefit from different types of control over changes to invariant subspaces. In particular, controlling the $2 \to \infty$ distance
\[
	\mathrm{dist}_{2 \to \infty}(V, \tilde{V}) :=
	\min_{W \in \Obb_{r}} \norm{V - \tilde{V} W}_{2 \to \infty}
	 = \min_{W\in \Obb_r} \max_{j \in [n]} \norm{V_{j, :} - \tilde{V}_{j, :} W},
\]
may be more applicable to settings where subspaces are interpreted row-wise (as
is the case in spectral clustering). Recent work on perturbation theory for the
$2 \to \infty$ norm~\cite{AFWZ17,CapeTangPriebe18,DamleSun19,EldBelWang18}
allows aspects our work to be extended to control incremental changes
to the subspace under this metric. However, the lack of convergence theory for
eigensolvers in this metric (beyond na\"ive bounds from norm equivalence)
prevent us from performing a full analysis in this setting.
\end{remark}

The eigensolvers sketched in~\cref{alg:subspace-iter,alg:block-krylov} perform a
prescribed number of iterations $k_{\max}$. However, in practice one usually
incorporates robust numerical termination criteria which often results in fewer
iterations being required. We omit such criteria here to simplify the presentation and defer
the details to the books by Parlett~\cite{Parlett98} and Saad~\cite{Saad11}.

\paragraph{Subspace iteration}
Subspace iteration, also known as simultaneous iteration, is a well-studied
generalization of the celebrated \textit{power method} for computing the leading
or trailing $r$ eigenvalues and eigenvectors of a symmetric matrix $A$. Starting
from an initial guess $V_0$, the algorithm computes orthonormal bases for terms
of the form $A^q V_0$, usually handled in practice via the QR factorization.
In the special case of $r = 1$, the intuition is that higher powers of $A$
amplify the component of the vector corresponding to the dominant eigenvector
(as long as the initial guess is not perpendicular to that eigenvector).
Since the estimates of subspace iteration converge to the eigenspace corresponding
to eigenvalues of largest \textit{magnitude} (instead of algebraically largest),
we assume it will be used with an appropriate shift.

\Cref{alg:subspace-iter} shows a full version of subspace iteration
incorporating the Rayleigh-Ritz procedure. In the algorithm, $r$
stands for the desired number of eigenpairs sought, but we may opt for running
it with a larger subspace of size $\ell > r$ for reasons that will be clarified later.
\begin{algorithm}[ht]
    \caption{Subspace iteration}
    \begin{algorithmic}[1]
        \State \textbf{Input}: matrix $A$, dimension $\ell$, target dimension $r$,
        initial estimate $V_0 \in \Rbb^{n \times \ell}$, $k_{\mathrm{max}}$.
        \For{$k = 1, \dots, k_{\max}$}
			\State $V_k, R_k := \mathtt{qr}(A V_{k-1})$
            \State $Q_k, \Lambda_k = \mathtt{spectral}({V_k}^\top A V_k)$
            \State $V_{k} := V_k \cdot Q_k$
        \EndFor
		\State \Return $(V_k)_{:, 1:r}, \diag(\Lambda_k)_{1:r}$
    \end{algorithmic}
    \label{alg:subspace-iter}
\end{algorithm}
There are a number of modifications (e.g., ``locking'' converged eigenpairs, shifts)
that one can incorporate to improve subspace iteration in practice~\cite{Saad11}.

The following theorem characterizes the convergence of
subspace iteration when $\ell = r$.
\begin{theorem}[Adapted from Theorem 8.2.2 in~\cite{GVL13}]
    \label{thm:subspace-iteration-convergence}
    Consider the spectral decomposition of $A$ as in \cref{eq:schur-decomp},
    assume that $\abs{\lambda_r} > \abs{\lambda_{r+1}}$, and let
    \(
        d_0 := \norm{V V^\top - V_0 V_0^\top}_2, \;
        \Delta_0 := \frac{d_0}{\sqrt{1 - d_0^2}},
    \)
    where $V_0 \in \Obb_{n, r}$.
    Then \cref{alg:subspace-iter} initialized with $V_0$ and $\ell = r$
    produces an orthogonal matrix $V_k$ such that
    \begin{equation}
        k_{\max} \geq \frac{\log(\Delta_0 / \varepsilon)}{\log\left(
            \abs{\frac{\lambda_r}{\lambda_{r+1}}}\right)}
        \Rightarrow \norm{V V^\top - V_k V_k^\top}_2 \leq \varepsilon.
        \label{eq:subit-iters}
    \end{equation}
\end{theorem}
The rate of convergence of individual eigenpairs may be faster than
the rate implied by \cref{thm:subspace-iteration-convergence}; for example,
the rate of convergence of the $i$\tsup{th} eigenvalue estimate is asymptotically
$\cO(\abs{\lambda_i / \lambda_{r+1}})$~\cite{Stewart1969}.

The amount of work suggested by \cref{thm:subspace-iteration-convergence}
crucially relies on two quantities: the eigenvalue ratio $\rho := \frac{\lambda_r}{\lambda_{r+1}}$,
which controls the convergence rate, as well as the initial subspace distance $d_0$.
In this paper, we use the fact that even if $\rho$ is moderately close to $1$, a small
initial distance can result in fewer iterations.

\paragraph{Block Krylov method}
Subspace iteration, while practical, discards a lot of useful information by
only using the last computed power $A^q V_{0}$, instead of
using the full (block) \textit{Krylov subspace} $\mathrm{span}(V_0, A V_0, \dots, A^{q}V_0)$.
This is precisely the motivation for (block) Krylov methods. The
well-known Lanczos method is part of this family when $V_0 = v_0$, a single
vector. In practice, the Lanczos method often exhibits superlinear convergence~\cite{Parlett98,Ruhe79},
making it the method of choice for eigenvalue computation.
However, its performance deteriorates in the presence of repeated
or \textit{clustered} (i.e., very close to each other in magnitude)
eigenvalues. It is well-known that the single-vector Lanczos method cannot find
multiple eigenvalues without deflation~\cite{Saad11} and converges slowly for
clustered eigenvalues. On the other hand, the block Lanczos method can get
around this issue, provided the block size $b$ is set appropriately (i.e., $b >
r$). See Parlett~\cite{Parlett98} for a discussion on the tradeoffs of block
sizes.

\Cref{alg:block-krylov} presents a Block Krylov method, albeit an impractical variant as it forms the
full Krylov matrix. While far from a practical implementation,
this form is helpful for stating and interpreting results under the assumption that we are operating in exact
arithmetic.
\begin{algorithm}[ht]
    \caption{Block Krylov method}
    \begin{algorithmic}[1]
        \State \textbf{Input}: $A \in \Rbb^{n \times n}$, dimension $r$,
        block size $b \geq r$, initial matrix $V_0 \in \Rbb^{n \times
        b}$, power $k_{\max}$.
        \State Form the block Krylov matrix \label{eq:krylov-matrix}
		\[
            K = \bmx{
                A V_0 & (A A^\top) A V_0 & \dots &
                (A A^\top)^{k_{\max}} A V_0
            }
		\]
        \State $Q, R = \mathtt{qr}(K)$
        \State $V_k, \Lambda_k = \mathtt{spectral}(Q^\top A Q)$
        \State \Return $(Q V_k)_{:, 1:r}, \diag(\Lambda_k)_{1:r}$
    \end{algorithmic}
    \label{alg:block-krylov}
\end{algorithm}

The block Lanczos algorithm first appeared to address the
shortcomings of the single-vector Lanczos iteration~\cite{CulDon74}. Saad
provides a comprehensive convergence analysis~\cite{Saad80}, while more recent work
presents sharp bounds and addresses convergence to clusters of eigenvalues~\cite{LiZhang15}.
For our purposes, we use the analysis found in Wang et al.~\cite{WZZ15}, the results
of which are stated in a form easily comparable with the corresponding rate
for~\cref{alg:subspace-iter}.

\begin{theorem}[Adapted from Theorem 6.6 in~\cite{WZZ15}]
    \label{thm:bk-convergence}
    Consider the Schur decomposition of $A$ as in \cref{eq:schur-decomp},
    assume that $\lambda_r > \lambda_{r+1}$, and let
    \(
    	d_0 := \norm{V V^\top - V_0 V_0^\top}_2,
    	\; \Delta_0 := \frac{d_0}{\sqrt{1 - d_0^2}}
    \), where $V_0 \in \Obb_{n, r}$. Then
    \cref{alg:block-krylov}
    initialized with $V_0$ produces an orthogonal matrix $V_k$ such that
    \begin{equation}
        k_{\max} \geq \frac{1 + \log_2(\Delta_0 / \varepsilon)}{
        \sqrt{\abs{\frac{\lambda_r}{\lambda_{r+1}}} - 1}}
        \Rightarrow
        \norm{V V^\top - V_k {V_k}^\top}_2 \leq \varepsilon.
        \label{eq:bk-iters}
    \end{equation}
\end{theorem}
Comparing with the rate of~\cref{thm:subspace-iteration-convergence}, the block
Lanczos method is clearly superior to subspace iteration: in the challenging
regime where $\lambda_r \approx \lambda_{r+1}$,
$\sqrt{\frac{\lambda_r}{\lambda_{r+1}} - 1} \gg
\log\left(\frac{\lambda_r}{\lambda_{r+1}}\right)$. As
with~\cref{alg:subspace-iter}, it is important to keep in mind that starting
from an estimate $V_0$ with $d_0 \ll 1$ can result in a significant reduction
of the required number of iterations.

\section{Subspace estimation with incremental updates} \label{sec:time-subspace-estimation}

In this section, we develop our incremental approach to spectral estimation,
motivated by data that changes over time. As stated before, we are interested in
maintaining an invariant subspace across ``small'' modifications to the data.
More specifically, we have a matrix $A \in \Rbb^{n \times n}$ and wish to
dynamically update an estimate of the \textit{leading} invariant subspace of $A$
corresponding to the eigenvalues of largest magnitude. This is not without loss
of generality; for example, modifying subspace iteration for finding the
smallest eigenpairs involves solving a linear system. However, some applications
stated in terms of finding the smallest eigenpairs can be reformulated to fit
our setting (see \Cref{sec:graph-spectral-embeddings}). A common example of such
a problem is spectral clustering for community detection in graphs, which we
rely on for intuition and constructing various numerical examples in the
experimental section. However, our approach generalizes beyond the graph
setting.

We formalize the incremental updates to our matrix as follows:
\begin{assumption} \label{asm:updates}
    Let $A_{(0)} := A \in \Rbb^{n \times n}$. At each time step $t \in \set{1,
    \dots, T}$, we observe updates of the form $A_{(t)} := A_{(t-1)} + E_{(t)}$,
    where $E_{(t)}$ is random, sparse, or low-rank. We denote
    the eigenvectors and eigenvalues of $A_{(t)}$ by $V_{(t)}$ and $\Lambda_{(t)}$.
\end{assumption}
An additional assumption that we impose on our matrix $A$ is that there is
sufficient decay on the lower end of the spectrum \textbf{outside} of the
subspace of interest, stated as \cref{asm:eigval-assumption}.
\begin{assumption} \label{asm:eigval-assumption}
    Let $\lambda_1 \geq \lambda_2 \geq \dots \geq \lambda_n$
    be the ordered eigenvalues of $A$. There is a constant $\gamma > 0$ and
    an integer $r_0 \ll n$ such that $\forall r \leq r_0$, there exists a
    $p \ll n$ satisfying
    \begin{equation}
        \abs{\lambda_r} \geq (1 + \gamma) \abs{\lambda_p}.
        \label{eq:eigval-assumption}
    \end{equation}
\end{assumption}
Intuitively, \cref{asm:eigval-assumption} implies a sort of
``uniform'' decay on the eigenvalues of $A$, since $\gamma$ does not depend on
the particular choice of $r \leq r_0$.
Indeed, if we allowed $\gamma = 0$, then \cref{eq:eigval-assumption} would be
trivially satisfied. The challenging regime is when $\gamma \ll 1$, as hinted
by standard convergence results on iterative methods for eigenvalue
problems~\cite{Demmel97,GVL13,Saad11}.

This technical assumption allows us to derive a priori upper bounds on the
required steps of our eigensolvers discussed above. In all the experiments
presented in \Cref{sec:applications}, \cref{asm:eigval-assumption} was verified
to hold with $\gamma > 0.1$ and $p \leq r + 10$.

\subsection{Outline of our approach}
We now sketch the high-level idea behind our approach. Naturally,
not all modifications to the underlying data have the same effect on the
invariant subspace of interest. For instance, ``small'' (with respect to some matrix norm) modifications should not
significantly change the leading eigenvalues and eigenvectors. Matrix
perturbation theory provides us with the analytical tools
to measure the worst-case behavior of these modifications~\cite{StewSun90}. Intuitively, if
this behavior is sufficiently ``bounded'', seeding our iterative method at
the previously obtained estimate should reduce the computation
required to obtain the next estimate. This is the crux of our approach.

An indispensable tool for perturbation bounds is the Davis--Kahan
theorem~\cite{DavKah70,Li98}. In order to state it formally, we first recall the
definition of \textit{principal angles} between subspaces~\cite{StewSun90}.
\begin{definition} \label{defn:principal-angles}
	Consider matrices $V, W\in \Obb_{n, k}$ and their corresponding column
	spaces $\cV, \cW$. Denote $\sigma_i := \sigma_i(V^\top W), \; i = 1, \dots, k$.
	The principal angles between $\cV$ and $\cW$ are the diagonal elements of the matrix
	\[
		\Theta(V, W) := \mathrm{diag}(\cos^{-1} \sigma_1, \dots, \cos^{-1} \sigma_k).
	\]
	Moreover, the $\sin$-distance between $\cV, \cW$ for any unitarily invariant norm $\norm{\cdot}$  is
	\begin{equation}
		\norm{\sin \Theta(V, W)} := \norm{\mathrm{diag}\left(
			\sin (\cos^{-1} \sigma_1), \dots, \sin (\cos^{-1} \sigma_k)
		\right)}.
		\label{eq:sin-distance}
	\end{equation}
\end{definition}

With \cref{defn:principal-angles} in hand, we can state the
Davis--Kahan \cref{thm:dav-kah}. While it holds for any unitarily
invariant norm, we are only interested in the case where $\norm{\cdot} =
\norm{\cdot}_2$.

\begin{theorem}[Davis--Kahan] \label{thm:dav-kah}
    Suppose $A, \hat{A}$ are symmetric and let $V, \hat{V} \in \Rbb^{n \times
    r}$ be two invariant subspaces containing eigenvectors of $A$ and $\hat{A}$,
    respectively. Let $\Lambda(V)$ denote the set of eigenvalues of $A$
    corresponding
    to $V$ and $\Lambda(\hat{V}_{\perp})$ the set of eigenvalues of $\hat{A}$
    corresponding to $\hat{V}_{\perp}$. If there exists an interval $[\alpha,
    \beta]$ and $\delta > 0$ such that $\Lambda(V) \in [\alpha, \beta]$ and
    $\Lambda(\hat{V}_{\perp}) \in (-\infty, \alpha - \delta) \cup
    (\beta + \delta, +\infty)$, then
    \begin{equation}
        \norm{\sin \Theta(V, \hat{V})}_2
        \leq \frac{\norm{(\hat{A} - A) V}_2}{\delta}.
        \label{eq:davkah}
    \end{equation}
    Moreover, if $\norm{A - \hat{A}}_2 < \frac{\delta_r}{2}$, where
    $\delta_r = \min\set{\abs{\mu - \lambda} \mmid \mu \in \Lambda(V_{\perp}), \lambda \in \Lambda(V)}$,
    then
    \begin{equation}
        \norm{\sin \Theta(V, \hat{V})}_2
        \leq \frac{2 \norm{(\hat{A} - A) V}_2}{\delta_r}.
        \label{eq:davkah-asm}
    \end{equation}
\end{theorem}
\Cref{eq:davkah-asm} above is straightforward to derive (e.g., see
the proof of~\cite[Corollary 2.1]{FWZZ18}). We employ \cref{thm:dav-kah}
to obtain an a priori bound on the distance between previously available
subspace estimates and the invariant subspace of the updated matrix.
\Cref{alg:high-level} outlines our high-level strategy.
\begin{algorithm}[ht]
	\caption{High-level algorithm for subspace updates}
    \begin{algorithmic}[1]
        \State \textbf{Input}: Matrix $A := A_{(0)}$, initial subspace size $r$,
        high-order $q \geq 1$, threshold $\varepsilon$.
        \State $\triangleright$ compute estimates $\hat{V}_{(0)} \in \Rbb^{n \times r}$ and
        $\hat{\lambda}_{r}, \dots, \hat{\lambda}_{r+q}$.
        \label{op:init-estimates}
        \For{$t = 1, \dots, T$}
            \State $A_{(t)} := A_{(t-1)} + E_{(t)}$
            \State $\bullet$ estimate $\norm{\sin \Theta(V_{(t)}, V_{(t-1)})}_2 \leq d_t$
            \label{op:davkah-proxy}
            \Comment{\Cref{sec:dav-kah-implem}}
            \If{$d_t > \varepsilon$}
                \State $\bullet$ update $\hat{V}_{(t)}, \hat{\lambda}_r$ to tolerance
                $\varepsilon$
                \label{op:iterative-refinement}
                \Comment{\Cref{sec:new-subspace}}
                \State $\bullet$ update high order eigenvalues $\hat{\lambda}_{r+1}, \dots, \hat{\lambda}_{r+q}$ to tolerance $\varepsilon$
                \label{op:high-order-update}
                \Comment{\Cref{sec:high-order}}
                \State $\bullet$ update subspace size $r$
                \Comment{\Cref{sec:size-updates}}
            \EndIf
        \EndFor
    \end{algorithmic}
    \label{alg:high-level}
\end{algorithm} \\

In \Cref{alg:high-level}, the constant $q \geq 1$  is an
``oversampling'' factor which is utilized as a heuristic for updating the size
of the subspace on the fly as well as for estimating the convergence rate
$\lvert\frac{\lambda_r}{\lambda_{r+1}}\rvert$ for our choice of eigenvalue algorithm.
Step~\ref{op:iterative-refinement} updates the subspace estimate while
step~\ref{op:high-order-update} updates the higher order eigenvalue estimates.
These updates are triggered when our proxy for \cref{eq:davkah-asm}
is above a specified threshold $\varepsilon$.

\begin{remark}[Extension for singular subspaces]
Our method is easily adapted for computing the leading singular subspace of a
rectangular matrix $A \in \Rbb^{m \times n}$. Via the standard
\textit{dilation trick}~\cite{GVL13,StewSun90}, we can form
\[
    S = \pmx{0 & A \\ A^\top & 0} \Rightarrow
    \lambda_i(S) = \pm \sigma_i(A),
\]
and given the (thin) SVD of $A$, $A = U \Sigma V^\top$, it is easy to verify
that the $i$\tsup{th} eigenvector of $S$ is the concatenation of the $i$\tsup{th} left and
right singular vectors, $u_i$ and $v_i$. Alternatively, if one is interested
only in the left or right singular subspace, it is also possible to run
\cref{alg:high-level} on $A^\top A$ or $A A^\top$ implicitly,
without ever forming the complete matrix\textemdash at the expense of standard
conditioning consequences.
\end{remark}

\subsection{Detailed description}
Below, we describe how to efficiently implement each of the aforementioned
steps, addressing the fact that our estimates at each step are inexact. We
discuss $d_t$ last, focusing on the updates in lines~\ref{op:iterative-refinement}
and~\ref{op:high-order-update} of \cref{alg:high-level} first.

\subsubsection{Computing the new subspace} \label{sec:new-subspace}
Let us first assume that $d_t$ from step~\ref{op:davkah-proxy} has already
been computed. Given the previous subspace estimate $\hat{V}_{(t-1)}
\in \Rbb^{n \times r}$, we can seed our eigensolver (\cref{alg:subspace-iter} or \cref{alg:block-krylov})
with $\hat{V}_{(t-1)}$. Assume that
\[
	\abs{\hat{\lambda}_r(A_{(t-1)}) - \lambda_r(A_{(t-1)})} \leq \varepsilon_1, \;
	\abs{\hat{\lambda}_{r+1}(A_{(t-1)}) - \lambda_{r+1}(A_{(t-1)})} \leq \varepsilon_2.
\]
Applying Weyl's inequality (\cref{lemma:weyl-ineq}), the rate controlling convergence will be at least
\begin{equation}
	\frac{\lambda_r(A_{(t)})}{\lambda_{r+1}(A_{(t)})}
	\geq
    \frac{\lambda_r(A_{(t-1)}) - \norm{E_{(t)}}_2}{\lambda_{r+1}(A_{(t-1)}) +
    \norm{E_{(t)}}_2}
	\geq \frac{\hat{\lambda}_r - \norm{E_{(t)}}_2 - \varepsilon_1}{
	\hat{\lambda}_{r+1} + \norm{E_{(t)}}_2 + \varepsilon_2} =: \rho_t
	\label{eq:convergence-rate},
\end{equation}
and \cref{thm:subspace-iteration-convergence,thm:bk-convergence} say that we have to set
(with $\Delta_t := d_t / \sqrt{1 - d_t^2}$):
\begin{equation}
	k_{\max} \geq
       	\begin{cases}
		\frac{\log\left(\Delta_t / \varepsilon
		\right)}{\log(\rho_t)} &
			\text{ for \Cref{alg:subspace-iter}} \\[2mm]
	        \frac{1 + \log_2\left(\Delta_t /
	\varepsilon\right)}{\sqrt{\rho_t - 1}} &
			\text{ for \Cref{alg:block-krylov}}
        \end{cases}
	\label{eq:k-max-range}
\end{equation}
in order to guarantee
$\mathrm{dist}(V_{(t)}, \hat{V}_{(t)}) \leq \varepsilon$.

In order to set $k_{\max}$, we need to know the error bounds $\varepsilon_1$ and
$\varepsilon_2$ on the accuracy of the previously computed eigenvalues
$\hat{\lambda}_r$ and $\hat{\lambda}_{r+1}$. For $\hat{\lambda}_r$, we can
apply \cref{thm:ritz-val-distance} for the estimate $\hat{V}_{(t-1)}$ which
gives us
\begin{equation}
	\abs{\hat{\lambda}_r(A_{(t-1)}) - \lambda_r(A_{(t-1)})}
	\leq \rho(A_{(t-1)}) \norm{\sin\Theta(\hat{V}_{(t-1)}, V_{(t-1)})}_2^2 \leq
	\rho(A_{(t-1)}) \varepsilon^2.
	\label{eq:lr-error}
\end{equation}
For the higher order eigenvalue estimate $\hat{\lambda}_{r+1}$, we can deduce a
similar worst-case bound (see \Cref{sec:high-order} for details):
\begin{equation}
	\abs{\hat{\lambda}_{r+1}(A_{(t-1)}) - \lambda_{r+1}(A_{(t-1)})}
	\leq 2\rho(A_{(t-1)}) \varepsilon^2.
	\label{eq:lr+-error}
\end{equation}
In \cref{eq:lr+-error}, the error factor is due to the combined inexactness of
the ``deflated'' matrix and the estimate returned from the iterative method of
choice.

\subsubsection{Computing the higher order eigenvalues} \label{sec:high-order}
The best way for updating the higher order eigenvalue estimates
$\hat{\lambda}_{r+1}, \dots$, $\hat{\lambda}_{r+q}$ is not obvious. One
possible approach is to ``augment'' the tracked subspace updated at
Line~\ref{op:iterative-refinement} to have dimension $r + q$. However, because
iterative eigenvalue algorithms use some orthogonalization scheme (e.g., QR
factorization) that scales quadratically in the subspace
dimension~\cite{GVL13}, we would have incurred a cost proportional to (at least)
$(r + q)^2$.

A more careful approach is a 2-phase algorithm, where we lock the
$r$-dimensional estimate $\mxhat{V}{t}$ and apply an iterative method to the
``deflated'' matrix
$(I - \mxhat{V}{t} \mxhat{V}{t}^\top) A_{(t)}(I - \mxhat{V}{t} \mxhat{V}{t}^\top)$,
with potentially substantially smaller time and memory costs than would be required
for maintaining an $(r+q)$-dimensional subspace. Even though it is evident that
deflating $A_{(t)}$ this way must incur some loss of accuracy,
by \cref{lemma:more-iterative-eigvals} this loss is negligible when
$\varepsilon \leq \rho(A_{(t)})^{-1}$.
\begin{lemma} \label{lemma:more-iterative-eigvals}
	Consider a matrix $A = \bmx{V & V_{\perp}} \bmx{\Lambda & 0 \\ 0 &
		\Lambda_{\perp}} \bmx{V & V_{\perp}}^\top$ and an estimate $\Vhat \in
	\Obb_{n, r}$ such that $\norm{V V^\top - \Vhat {\Vhat}^\top}_2 \leq \varepsilon$.
        Then
	\begin{align}
	\abs{\lambda_{r+i}(A) - \lambda_i\left( (I - \Vhat {\Vhat}^\top) A
		(I - \Vhat {\Vhat}^\top) \right)} &\leq \rho(A) \varepsilon^2
	\end{align}
\end{lemma}
\begin{proof}
	If $V$ contains a basis for the leading $r$-dimensional eigenspace, then
	it is obvious that $I - V V^\top$ is the projector to the trailing
	$(n - r)$-dimensional eigenspace, which corresponds to a contiguous set of
	eigenvalues of $A$. By elementary arguments, we know that
	\[
		\lambda_{r+i}(A) = \lambda_i\big(\underbrace{
        (I - V V^\top) A (I - V V^\top)}_{ =: \tilde{A}} \big).
	\]
	Therefore, by \cref{thm:ritz-val-distance},
	\begin{align*}
	\max_{i} \abs{\lambda_i\left(\tilde{A}\right)
		- \lambda_i\left((I - \Vhat {\Vhat}^\top) A (I - \Vhat {\Vhat}^\top)
		\right)}
	\leq \rho(A) \norm{V V^\top - \Vhat {\Vhat}^\top}_2^2
	\leq \rho(A) \varepsilon^2,
	\end{align*}
	since $(I - \Vhat {\Vhat}^\top)$ is also a projector to an $(n-r)$-dimensional
	subspace.
\end{proof}
We conclude the following accuracy estimate (dropping the $(t)$ subscript for brevity):
\begin{align}
	\begin{aligned}
	\abs{\lambda_{r+i}(A) - \hat{\lambda}_{r+i}(A)} &\leq
	\abs{\lambda_{r+i}(A) - \lambda_i\left(
		(I - \hat{V} \hat{V}^\top)A(I - \hat{V} \hat{V}^\top)\right)} \\
	&\quad
		+ \abs{\lambda_i\left(
		(I - \hat{V} \hat{V}^\top)A(I - \hat{V} \hat{V}^\top)\right)
		- \hat{\lambda}_i\left(
		(I - \hat{V} \hat{V}^\top)A(I - \hat{V} \hat{V}^\top)\right)} \\
	& \leq 2 \rho(A) \varepsilon^2
	\end{aligned},
	\label{eq:lr+-estimate}
\end{align}
with the first $\rho(A) \varepsilon^2$ factor coming from an application
of \cref{lemma:more-iterative-eigvals} and the second such factor coming from an
application of \cref{thm:ritz-val-distance} when applying our iterative method
with accuracy parameter $\varepsilon$ to the deflated matrix. This provides the
inequality in \cref{eq:lr+-error}.

When computing higher order eigenvalues, we might not have a previously
maintained estimate of the corresponding subspace. In this situation,
it is common to pick a random Gaussian matrix as the seed matrix. The
two following Propositions provide guarantees for
\cref{alg:subspace-iter,alg:block-krylov} under this initialization scheme.
\begin{proposition}[Corollary of Theorem 5.8 in~\cite{Gu15}]
	\label{prop:rand-subspace-guarantee}
	Let $A = V \Lambda V^\top$ be a symmetric matrix, and let
	$0 \leq p \leq \ell - k$. For a given $\delta \in (0, 1)$, define
	\[
		C_{\delta} := \frac{e \sqrt{\ell}}{p + 1}
		\left(\frac{2}{\delta}\right)^{\frac{1}{p+1}}
		\left[
		\sqrt{n - \ell + p} + \sqrt{\ell} + \sqrt{2
			\log\left(\frac{2}{\delta}\right)}
		\right].
	\]
	Then, with probability at least $1 - \delta$, for $j = 1, \dots, k$,
	the Ritz values $\set{\mu_j}$ returned by \cref{alg:subspace-iter}
	initialized with a random Gaussian matrix $(V_0)_{ij} \sim \cN(0, 1)$
	satisfy
	\[
		\lambda_j \geq \mu_j \geq
		\frac{\lambda_j}{\sqrt{1 + C_{\delta}^2
				\left(\frac{\lambda_{\ell - p + 1}}{\lambda_j}\right)^{4k_{\max}
					+ 2}}},
		\quad \forall j = 1, \dots, r.
	\]
\end{proposition}

\begin{proposition}[Theorem III.4 in~\cite{YuanGuLi18}]
	\label{prop:rand-block-krylov}
	Let $A = V \Lambda V^\top$ be a symmetric matrix, and let $b = r + p$
	in \cref{alg:block-krylov} initialized with $(V_0)_{ij} \sim \cN(0, 1)$.
	Then, for $j = 1, \dots, r$, the Ritz values $\set{\mu_j}$ returned satisfy
	\[
		\lambda_j \geq  \mu_j \geq
		\frac{\lambda_j}{\sqrt{1 + C^2 T_{2k_{\max}+1}^{-2} \left(1 +  2
				\frac{\lambda_j - \lambda_{j + p +1}}{\lambda_{j + p +1}} \right)}},
	\]
	where $T_{p}$ is the $p$\tsup{th} degree Chebyshev polynomial and $C$ is a
	constant depending on the initialization matrix.
\end{proposition}

Let us incorporate \cref{asm:eigval-assumption} into
\cref{prop:rand-subspace-guarantee,prop:rand-block-krylov}. For the latter, it is known that Chebyshev
polynomials satisfy $T_{d}(1 + \alpha) \geq 2^{d \sqrt{\alpha}}$ when
$\alpha \ll 1$ (see~\cite[Lemma 5]{MusMus15}). Therefore, to achieve $\varepsilon$-accuracy under
Gaussian initialization, we can set
\begin{equation}
	k_{\max} \geq
	\begin{cases}
		\frac{\log(C_1 / \varepsilon)}{\log(1 + \gamma)} & \text{ for \cref{alg:subspace-iter}} \\[2mm]
		\frac{\log(C_2 / \varepsilon)}{\sqrt{\gamma}} & \text{ for \cref{alg:block-krylov}},
       \end{cases}
	\label{eq:kmax-gaussian-init}
\end{equation}
where $C_1, C_2$ are the constants from \cref{prop:rand-subspace-guarantee,prop:rand-block-krylov}.

\subsubsection{A proxy for the Davis--Kahan bound} \label{sec:dav-kah-implem}
We still need a good bound $d_t$ on the subspace distance in \cref{eq:davkah-asm}.
Evaluating $\norm{(A_{(t)} - A_{(t-1)})V}_2$ exactly is not possible, since
we only maintain approximations of $V$. A standard majorizer, often too
loose in practice, is $\norm{A_{(t)} - A_{(t-1)}}_2$. We can get
a better upper bound by leveraging a priori information about the
perturbation, such as sparsity or randomness\textemdash two such instances are handled
in \Cref{sec:applications}.

We can still get around this issue in the absence of prior information about
$E_{(t)} = A_{(t)} - A_{(t-1)}$. Let us first assume that we have computed a
$\varepsilon$-close subspace estimate $\mxhat{V}{t}$ (i.e.,
$\mathrm{dist}(V_{(t)}, \mxhat{V}{t}) \leq \varepsilon$
).
Given this information, we can compute an upper bound to the subspace distance
at time $t+1$ in terms of $\mxhat{V}{t}$ and the perturbation $E_{(t)}$, as \cref{lemma:approximate-davis-kahan} shows. Again, we drop the subscripts
for simplicity.
\begin{lemma} \label{lemma:approximate-davis-kahan}
    Let $E$ be a perturbation to a matrix $A = \bmx{V & V_{\perp}}
    \Lambda \bmx{V & V_{\perp}}^\top$ and $\Vhat$ be
    a subspace estimate such that
    \(
        \norm{V V^\top - \Vhat {\Vhat}^\top}_2 \leq \varepsilon.
	\)
    Then
    \begin{align}
        \norm{E V}_2 &\leq
        \sqrt{ \norm{E \Vhat}^2_2 + \varepsilon \norm{E}^2_2 }.
        \label{eq:approx-dav-kah-numerator}
    \end{align}
\end{lemma}
\begin{proof}
    Using the triangle inequality for the first inequality and our assumption
    on $V, \hat{V}$ for the second inequality, we obtain
    \begin{align*}
        \norm{E V}_2 &\leq \sqrt{\norm{E (VV^\top - \Vhat {\Vhat^\top})
        E^\top}_2
        + \norm{E \Vhat {\Vhat}^\top E^\top}_2}
        \leq \sqrt{\varepsilon \norm{E}^2_2 + \norm{E \Vhat}_2^2}.
    \end{align*}
\end{proof}

This provides us with a proxy for the numerator of \cref{eq:davkah}, with the
estimate for the eigenvalue gap following from the discussion
in \Cref{sec:new-subspace,sec:high-order}. From \cref{eq:lr-error,eq:lr+-error},
we know that the approximate eigengap
$\hat{\delta}_r := \hat{\lambda}_r - \hat{\lambda}_{r+1}$ satisfies
\begin{align*}
    \abs{\hat{\delta}_r - \delta_r} &=
    \abs{\lambda_r({V_0}^\top A V_0) - \lambda_1(
    (I - V_0 {V_0}^\top) A (I - V_0 {V_0}^\top))} \leq
    3\varepsilon^2 \rho(A).
\end{align*}

Putting everything together, we arrive at \cref{corl:dav-kah-proxy},
whose proof follows immediately after
applying \cref{lemma:approximate-davis-kahan} to bound the numerator on
\cref{thm:dav-kah}.
\begin{corollary} \label{corl:dav-kah-proxy}
    Given a sequence of updates satisfying \cref{asm:updates},
    \cref{alg:high-level} maintains an estimate $\Vhat_{(t)}$
    satisfying $\mathrm{dist}(\Vhat_{(t)}, V_{(t)}) \leq \varepsilon$.
    If $\abs{\lambda_{r}(A_{(t-1)}) - \lambda_{r+1}(A_{(t-1)})} <
    \frac{\norm{E_{(t)}}_{2}}{2}$, Line~\ref{op:iterative-refinement} of
    \cref{alg:high-level} takes at most $k_{\max}$ iterations, with
    \begin{align}
       	k_{\max} \geq
       	\begin{cases}
       		\frac{\log(\Delta_t / \varepsilon)}{\log(\rho_t)} &
       			\text{ for \cref{alg:subspace-iter}} \\[2mm]
       		\frac{1 + \log(\Delta_t / \varepsilon)}{\sqrt{\rho_t - 1}} &
       			\text{ for \cref{alg:block-krylov},}
       	\end{cases}
        \label{eq:kmax-range-corl}
    \end{align}
    where $\Delta_t := \frac{d_t}{\sqrt{1 - d_t^2}}$ and $\set{d_t, \rho_t}$ are given by
    \begin{align}
        \begin{aligned}
        d_t &:= \frac{2 \sqrt{\varepsilon \norm{E_{(t)}}^2_2 +
        \norm{E_{(t)}
                \Vhat_{(t-1)}}^2_2}}{
                \hat{\lambda}_{r}(A_{(t-1)}) -
                \hat{\lambda}_{r+1}(A_{(t-1)}) - 3\varepsilon^2
                \rho(A_{(t-1)})} \\
        \rho_t &:= \frac{\hat{\lambda}_{r}(A_{(t-1)}) - \norm{E_{(t)}}_2
        - \rho(A_{(t-1)}) \varepsilon^2}{
        \hat{\lambda}_{r+1}(A_{(t-1)}) + \norm{E_{(t)}}_2
        + 2 \rho(A_{(t-1)}) \varepsilon^2}.
        \end{aligned}
        \label{eq:rho_t-delta_t}
    \end{align}
    Moreover, under \cref{asm:eigval-assumption},
    Line~\ref{op:high-order-update} takes at most $k_{\max}'$
    iterations, with
    \begin{equation}
    	k_{\max}' \geq
    	\begin{cases}
    		\frac{\log(C_1 / \varepsilon)}{\log(1 + \gamma)} & \text{
    		for \cref{alg:subspace-iter}} \\[2mm]
    		\frac{\log(C_2 / \varepsilon)}{\sqrt{\gamma}} & \text{
    		for \cref{alg:block-krylov},}
    	\end{cases}
        \label{eq:kmax-hi}
    \end{equation}
    where $C_1, C_2$ are the constants
    in \cref{prop:rand-subspace-guarantee,prop:rand-block-krylov}.
\end{corollary}

\begin{remark}
	Though the factor $\rho(A)$ above seems to introduce a non-negligible loss,
	several applications of interest satisfy $\rho(A) = \cO(\varepsilon^{-1})$.
    For example, spectral clustering with the
	normalized adjacency matrix satisfies $\rho(A) \leq 2$. On the numerical
	side, if we can't assume that $\rho(A) = \cO(\varepsilon^{-1})$, it is
	sometimes possible to estimate $\rho(A_{(t)}) \leq
	\hat{\lambda}_1(\mxhat{A}{t-1}) + \norm{E_{(t)}}_2$
	using \cref{lemma:weyl-ineq}, as in the case of positive semi-definite
	matrices. This task is straightforward in all applications considered
	in \Cref{sec:applications}.

    If the assumption of the Davis-Kahan theorem fails or
    $d_t \geq 1$, it is straightforward to check that seeding with a
    Gaussian matrix results in the bound of \cref{eq:kmax-hi}
    for Line~\ref{op:iterative-refinement}, with $\rho_t$ instead of $\gamma$
    if $\rho_t > 1$.
\end{remark}

\subsubsection{Subspace size updates} \label{sec:size-updates}
As discussed in \Cref{sec:intro-temporal}, it is often the case
that the dimension $r$ of the subspace of interest is either not known a priori
or changing over time. In spectral clustering, a common heuristic to determine
the ``correct'' size of clusters is to compute a large extremal eigenspace of
the normalized adjacency matrix and observe the ratios of successive
eigenvalues~\cite{VonLux07}.

We can extend this criterion to enable tracking the dimension over time, as
approximating higher order eigenvalues given an $r$-dimensional invariant
subspace is feasible via the procedure described in
\cref{lemma:more-iterative-eigvals}. Assume that estimates
$\hat{\lambda}_1, \dots, \hat{\lambda}_r, \dots, \hat{\lambda}_{r + q}$ are
available. We propose setting the new candidate size $r_{c}$ as follows:
\begin{equation}
    r_{c} := \argmin_{2 \leq i \leq q - 1}
    \abs{\frac{\hat{\lambda}_{i+1}}{\hat{\lambda}_i}}.
    \label{eq:candidate-size}
\end{equation}
In other words, $r_c$ is set to the index of the smallest ratio of successive
eigenvalues, with $i \geq 2$ in \cref{eq:candidate-size} to exclude trivial
subspaces with just $1$ element. We may also adapt \cref{eq:candidate-size} to
take into account the recent history of eigenvalue ratios, setting $r_c \neq r$
only after $T$ rounds of $\hat{\lambda}_{r_c + 1} / \hat{\lambda}_{r_c}$
as the dominant eigenvalue ratio to account for short-lived fluctuations.

\section{Applications} \label{sec:applications}

In this section, we illustrate the effectiveness of our incremental subspace
estimation framework with two applications. The first application is tracking
the leading subspace of the normalized adjacency of a graph evolving in time,
which enables clustering the underlying graph incrementally. In this case, the
changes are \textit{sparse}, since only a small fraction of the edges of the
graph changes over time. The data matrices are also sparse, since the graphs
involved have small vertex degrees, so matrix-vector multiplication is
efficient. The second application is Principal Component Analysis (PCA), with
deterministic or random low-rank updates. In this case, the data matrices are
dense, but in some cases, we can use structure (such as Toeplitz/Hankel
structure) to speed up matrix-vector multiplication.  For both applications, we
outline domain-specific properties which simplify some of the involved
computations and can further improve the performance of our adaptive procedure.

Numerical experiments for \cref{alg:block-krylov} use
the locally optimal block preconditioned conjugate gradient method
(\texttt{LOBPCG})~\cite{Knyazev01}, which is the only available block Krylov
implementation in \texttt{Julia} that allows seeding the subspace.
\texttt{LOBPCG} also provides built-in support for orthogonalization
``constraints'' which are used in Line~\ref{op:high-order-update}
of \cref{alg:high-level} to implement matrix-vector multiplication with the
deflated matrix. We do not use a preconditioner when calling
\texttt{LOBPCG}. Finally, all of the code is freely available
from~\texttt{\url{https://github.com/VHarisop/inc-spectral-embeddings}}.

\subsection{Spectral embeddings in graphs} \label{sec:graph-spectral-embeddings}
A common problem in analyzing graph data is finding its clusters, communities,
or modules. A standard approach to this problem is \textit{spectral
clustering}~\cite{NgJorWei02,VonLux07}. This procedure applies the $k$-means
algorithm with an initial seed $V_0 \in \Rbb^{d \times r}$, where $V_0$ is the
invariant subspace of the normalized adjacency matrix of the graph corresponding
to the $r$ largest eigenvalues. There is a fundamental computational problem for
spectral clustering when the graph is changing over time: if the graph is large
enough, it simply becomes too expensive to compute the embedding from scratch
every time. However, we often observe individual ``small'' (low-rank or small
norm) and sparse perturbations to the graph, such as the addition of an
edge. Therefore, our incremental spectral estimation framework is well-suited
for this task.

For a bit more notation, we denote a general undirected graph by $\cG = (\cV, \cE)$ and let
$n := \abs{\cV}$ and $m := \abs{\cE}$. We write $A$ for the adjacency matrix of
$\cG$. The \textbf{degree} matrix $D = \mathrm{diag}(A\mathbf{1})$, where $\mathbf{1}$ is the
vector of all ones, is the diagonal matrix whose diagonal entries contain the
degrees of each node. The \textbf{symmetrically normalized adjacency} matrix is then
$\tilde{A} = D^{-1/2} A D^{-1/2}$. Again, basic spectral clustering computes a
principal subspace of $\tilde{A}$ and runs $k$-means on this spectral embedding
to cluster the nodes.

In typical data, the adjacency matrix is sparse, and in this case, it is
possible to deduce a priori bounds which can aid us in the application of the
Davis--Kahan theorem as well as in the estimation of the convergence rate of the
iterative methods. To formalize this, assume that the adjacency matrix $A$ is
perturbed by a symmetric matrix $E \in \set{0, \pm 1}^{n \times n}$ such that
the number of modifications at each vertex does not exceed the number of its
currently incident edges. In that case it is possible to bound the
norm of the perturbation to the normalized adjacency matrix, as in
\cref{prop:perturb-norm-sparse-updates}.

\begin{proposition} \label{prop:perturb-norm-sparse-updates}
	Suppose we observe a sequence of edge updates corresponding to a
	symmetric matrix $E \in \set{0, \pm 1}^{n \times n}$ such that
	$\frac{1}{d_i} \abs{\sum_{j=1}^n E_{ij}} \leq \alpha < 1, \; i = 1, \ldots,
	n$
	and $E_{ii} = 0$. Let $\tilde{A}_{\mathrm{new}}$ denote the updated
	normalized adjacency matrix. Then
	\begin{equation}
        \norm{\tilde{A}_{\mathrm{new}} - \tilde{A}}_2 \leq
        \alpha \cdot \left(1 + \alpha + \sqrt{\min\set{\kappa(D), \rank(E)}
        }(1 + \cO(\alpha^2))\right) + \left(\frac{\alpha
        (1 + \alpha)}{2}\right)^2,
        \label{eq:sparse-update-bound}
	\end{equation}
    where $\kappa(D) := \frac{\max_i d_i}{\min_j d_j}$.
\end{proposition}
The proof is in \Cref{sec:sparse-updates-proof}.
\cref{prop:perturb-norm-sparse-updates} provides theoretical justification for
the intuition that updating vertices with few neighbors tends to have more
severe effects on the spectrum; in contrast, if all affected vertices have
large neighborhoods, we expect that $\alpha \ll 1$.

\subsubsection{Stochastic block model}
In this experiment, we use synthetic data generated by the stochastic block
model~\cite{HolLasLein83}. The probabilistic model consists of a set of $n$ nodes and $k$
clusters or communities $\set{\cC_i}_{i=1}^k$, with each node belonging to
exactly one community $\cC_i$. Edges are observed according to the following
model:
\begin{equation}
    \prob{v_i \text{ connected to } v_j} = \begin{cases}
        p, & v_i, v_j \in \mathcal{C}_p \text{ for some } p \\
        q, & \text{otherwise}
    \end{cases}.
    \label{eq:sbm-edges}
\end{equation}
However, we only observe a graph and not the cluster idenitification.  There is
a rich literature on how spectral clustering methods can identify the latent
communities~\cite{McSherry-2001-spectral,Rohe-2011-spectral,Sussman-2012-spectral}.

\begin{figure}[tb]
    \centerline{\includegraphics[width=\textwidth]{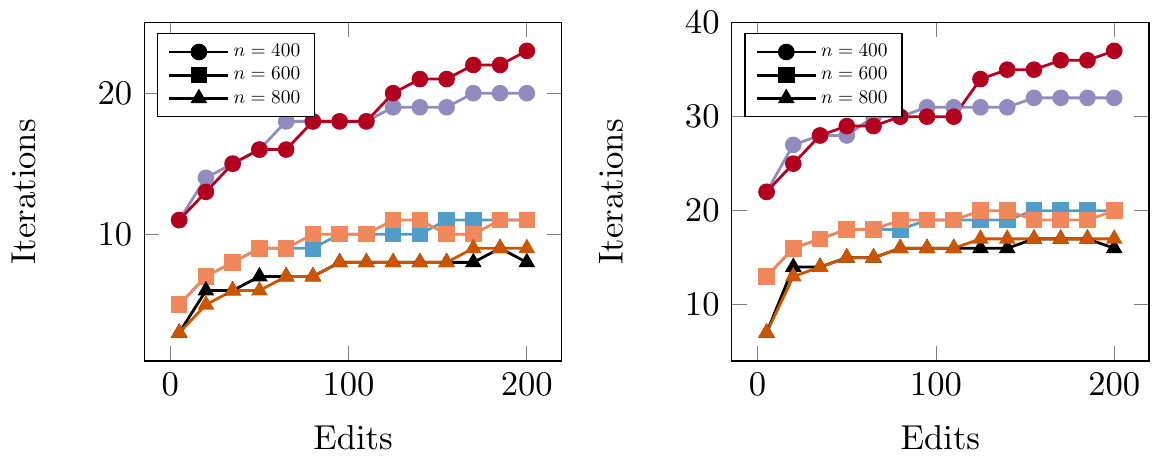}}
    \caption{Upper bound from~\eqref{eq:kmax-range-corl} on iterations
    performed for tracking $V \in \Rbb^{n \times 5}$ using
    \cref{alg:subspace-iter} without (cold colors) vs.\ with
    (warm colors) eigengap oracle. The resulting bounds exhibit negligible
    differences.
    Underlying graphs are SBMs with $p = 0.5, q = 0.1$. Left:
    $\varepsilon = 10^{-2}$. Right: $\varepsilon = 10^{-3}$.}
    \label{fig:SBM-iters-vs-oracle}
    \centerline{\includegraphics[width=\textwidth]{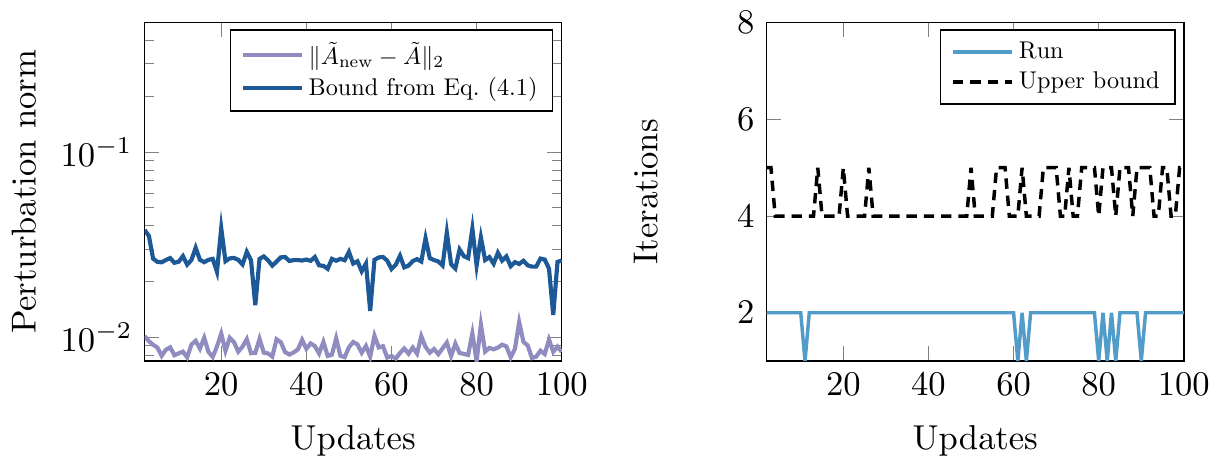}}
    \caption{Performance of \cref{alg:high-level} for tracking $V \in \Rbb^{n
    \times 5}$ using the bound from \cref{eq:sparse-update-bound} with
    $\alpha < 0.05$ (left) and resulting iterations (right). The perturbation
    bound is at most an order of magnitude off, with the resulting iteration
    bound just a constant multiple of those elapsed to
    reach accuracy $\varepsilon = 10^{-3}$. Underlying graphs are SBMs with
    $p = 0.5, q = 0.1, n = 1000$.}
    \label{fig:SBM-sparse-updates}
\end{figure}

\paragraph{Effect of eigenvalue precision}
We initially examine how the looseness of our eigenvalue estimates
affects (via the quantities $d_t, \rho_t$ in \cref{eq:rho_t-delta_t}) the bound
on the number of iterations. Our experimental methodology is as follows.
\begin{itemize}
\item We first sample two graphs $G_s$ and $G_t$ from the stochastic block model with
parameters $p = 0.5$, $q = 0.1$, with 5 clusters in $G_s$ and 6 clusters in $G_t$.
\item We then use a ``network interpolation'' method~\cite{ReeDamBen19} to
evolve $G_s$ towards $G_t$. At each step, we perform one edge addition or
deletion on $G_s$. This update decreases the graph edit distance
by $1$ with probability $h$, and increases it with probability $1 - h$.
\item Each update $E_{(t)}$ corresponds to 5 of the aforementioned edits, to be
used in single iteration of \cref{alg:high-level}. The eigensolver is run
for the full number of iterations prescribed by the upper bounds instead of
accuracy $\varepsilon$.
\end{itemize}
\Cref{fig:SBM-iters-vs-oracle} shows the number of iterations needed by
\cref{alg:high-level} as a function of the number of graph edits.
Warm colors correspond to using an oracle for the eigenvalues involved,
i.e., the real values of ${\lambda}_i$ (up to numerical accuracy) instead
of the estimates $\hat{\lambda}_i$, which correspond to the cold colors.
Interestingly, the difference is negligible, so we
should not expect our numerical bounds to degrade significantly with ``rough''
eigenvalue estimates except in challenging regimes where the ratio controlling
convergence is already very close to $1$.

All instances generated for this experiment interpolate between a graph with
five clusters and one with six clusters, so we expect the eigengap
$\lambda_{r+1} - \lambda_r$ to gradually decrease; therefore the
lowest-dimensional example ($n = 400$) should be the most challenging, which is
empirically verified. Smaller values of $\varepsilon$ also need fewer subspace
iterations, as expected. Finally, by comparing the prescribed accuracy
$\varepsilon$ with the final distance of the computed estimates,
$\mathrm{dist}(V_{(t)}, \mxhat{V}{t})$, we find that the latter is at least $2$
orders of magnitude smaller, a somewhat unsurprising outcome given the
pessimistic nature of most perturbation results used above.

\paragraph{Effect of sparse updates to spectrum}
We take a further step to assess the impact of sparse updates via the lens
provided by~\cref{prop:perturb-norm-sparse-updates}. More precisely, we generate
an SBM with $n = 1000$ nodes and $5$ communities, with $p = 0.5$, $q = 0.1$ as
previously. We then generate a sequence of $100$ random updates, each
corresponding to $50$ edge modifications. As clusters are sizeable and the
intra-cluster edge probability $p$ is moderately high, each of these updates
results in $\alpha < 0.05$ in the notation
of~\cref{prop:perturb-norm-sparse-updates}.

The performance of \cref{alg:high-level} using the worst-case estimate
of \cref{eq:sparse-update-bound} to majorize the Davis--Kahan bound is
illustrated in \cref{fig:SBM-sparse-updates}. The estimate is at most an order
of magnitude higher than the true perturbation norm. Moreover, the resulting
bound on iterations is off by a constant factor. Therefore, in situations where
we observe sparse updates to graphs with few isolated nodes, it is
possible to get additional speedups by replacing $d_t$
from \cref{eq:rho_t-delta_t} with a less complicated proxy.

\subsubsection{Time-evolving real-world graphs}
\label{sec:experiment-realgraphs}
We next test our algorithms on real-world graph datasets:
the \texttt{college-msg}~\cite{PanOpsCar09} dataset of private communications on
a Facebook-like college messaging platform, as well as subsets of
the \texttt{temporal-reddit-reply} dataset~\cite{LiuBenCha19} consisting of
replies between users on the public social media platform \texttt{reddit.com}.
All datasets are anonymized and contain timestamped edges representing the
interactions. We make edges undirected and remove any duplicates during
preprocessing. We isolate the largest connected component and work on the
induced subgraph for each dataset, leading to the statistics shown
in \cref{table:datasets}. For the \texttt{reddit-*} datasets, \texttt{small},
\texttt{medium} and \texttt{large} variants correspond to keeping the first
$N = \set{100000, 250000, 1000000}$ nodes of the raw data, respectively.

\begin{table}
    \centering
    \caption{Summary statistics of temporal graphs used in numerical experiments.}
    \begin{tabular}{r @{\quad} c c}
    \toprule
        \textbf{Dataset} & {\bf \# of Nodes} & {\bf \# of Edges} \\ \midrule
        \texttt{college-msg} & 1,893 & 13,834 \\
        \texttt{reddit-small} & 78,529 & 455,864 \\
        \texttt{reddit-medium} & 217,286 & 1,698,265 \\
        \texttt{reddit-large} & 757,015 & 5,487,069 \\ \bottomrule
    \end{tabular}
    \label{table:datasets}
\end{table}
As before, we isolate the subgraph $G = (\cV, \cE)$ corresponding to the largest
connected component of the static version of the graph, and start from a version
$G_s = (\cV_s, \cE_s)$ such that $\cV_s = \cV,\cE_s \subset \cE$, containing the
$\abs{\cE_s}$ edges with the earliest timestamps. Then, we add edges in the
order they were encountered in the original dataset.

We employ a regularized version of the normalized adjacency
matrix~\cite{ACBL13,JosYu16,QinRohe13,ZhangRohe18}. Specifically, given a
regularization parameter $\tau$, the regularized adjacency matrix is
$A_{\tau} := A + \frac{\tau}{n} \bm{1} \bm{1}^\top$.
We set $\tau$ equal to $1$ for the \texttt{college-msg} dataset and
equal to the average degree ($\tau = \frac{1}{n} \sum_{i=1}^n D_{ii}$),
for the \texttt{reddit-*} datasets.
This regularization can improve spectral clustering by addressing the adverse
effects of isolated or low-degree nodes.
In particular, the large cluster of the leading eigenvalue (corresponding to
connected components in the graph), can significantly degrade the performance
of both iterative methods.

\begin{figure}[tb]
    \centerline{\includegraphics[width=\textwidth]{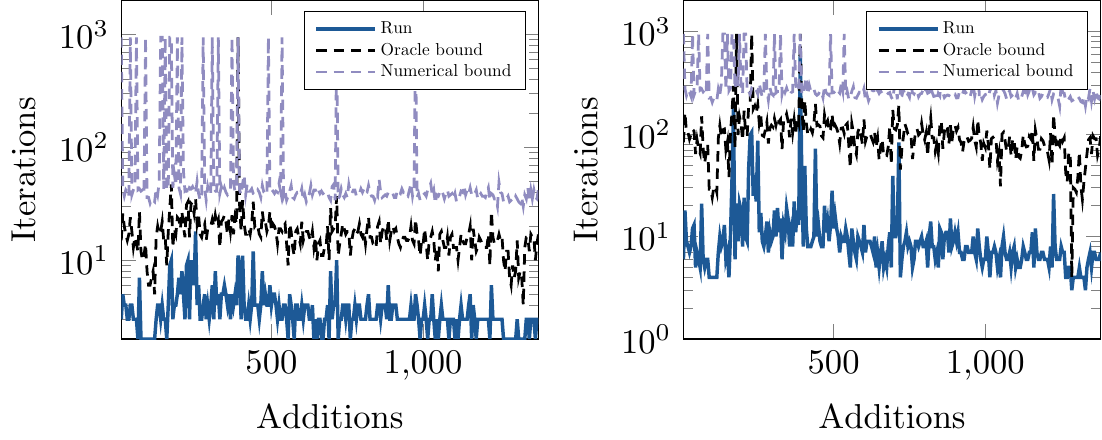}}
    \caption{Performance of \cref{alg:high-level} using block
    Krylov method (left) and subspace iteration (right) on the
    \texttt{college-msg} dataset, with accuracy $\varepsilon = 10^{-3}$. Both
    plots indicate that the bound of~\eqref{eq:kmax-range-corl} can be
    \textit{sharp}, even though it is loose in general.}
    \label{fig:dynamic-vs-oracle}
    \vspace{11pt}
    \centerline{\includegraphics[width=\textwidth]{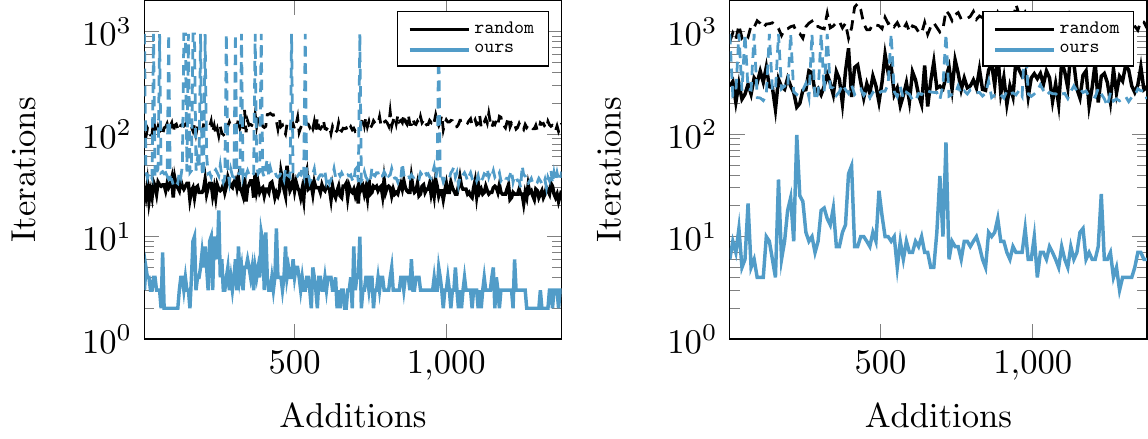}}
    \caption{Iterations using the Davis--Kahan bound
    and \cref{alg:subspace-iter} (left) vs.\ \cref{alg:block-krylov} (right),
    starting from $V_0$ (\texttt{ours}) vs.\ a random initial estimate
    (\texttt{random}) on the \texttt{college-msg} dataset ($\varepsilon = 10^{-3}$).
    Solid lines are iterations until numerical convergence,
    dashed lines are upper bounds. Naive seeding incurs orders of magnitude
    higher computational costs, which can even exceed the (conservative) upper
    bound of \cref{eq:kmax-range-corl}.}
    \label{fig:prediction-comparison}
\end{figure}

At each step, we use the previous estimate $\hat{V}$ of the leading subspace of
$A := A_{\tau}$ to initialize subspace iteration. Moreover, if we detect a
change in the dimension of the invariant subspace (using the criterion
in \Cref{sec:size-updates}), we follow the steps below to form
$\hat{V}_{\mathrm{new}}$:
\begin{itemize}
\item if $r_{\mathrm{new}} < r$, we simply drop the eigenvectors $v_i$
corresponding to the smallest Ritz values $\theta_i = \lambda_i(\hat{V}^\top
A \hat{V})$.
\item if $r_{\mathrm{new}} > r$, we generate $r_{\mathrm{new}} - r$ new vectors
and orthogonalize them against $\hat{V}$ using the (modified) Gram-Schmidt
process~\cite{GVL13}.
\end{itemize}

Initially, we experiment with the \texttt{college-msg} dataset, as it is
feasible to evaluate the true subspace distances and to compare with the
performance of the randomly initialized variant at each step.
\Cref{fig:dynamic-vs-oracle} shows the number of iterations run to
attain the desired accuracy (checked using numerical stopping criteria) of
$\varepsilon = 10^{-3}$, as well as the number of iterations $k_{\max}$
determined by \cref{eq:kmax-range-corl}, starting from $G_s$ and evolving
towards the static version of the graph.
As in the SBM case, we observe five edge modifications at a time.
The plot reveals that our bounds are essentially sharp, since there is more than
one occasion where the number of subspace iterations almost matches the
upper bound (appearing as spikes).

\newcommand{\drnd}{d_{\mathrm{rnd}}}
We also observe that our warm-starting methods
provides substantial performance gains compared to random initialization.
\Cref{fig:prediction-comparison} compares our pipeline
with randomly initialized subspace iteration. We set $d_t := \mathrm{dist}(V,
\hat{V})$, with $\hat{V}$ being the Q factor from the QR
decomposition of a random Gaussian matrix, and report the number of iterations
elapsed until achieving accuracy $\varepsilon$ starting from $\hat{V}$ as well
as the bounds prescribed
by~\cref{thm:subspace-iteration-convergence,thm:bk-convergence} using $\rho :=
\frac{\lambda_{r+1}(A)}{\lambda_r(A)}$ (computed using
\texttt{Arpack}). Our upper bound is usually 1--2 orders of
magnitude tighter than the bound using random initialization, while far fewer
iterations are needed to reach convergence when initializing $V_0$ with the
previous estimate.

\begin{figure}[tb]
	\centering
	    \includegraphics[width=0.32\linewidth]{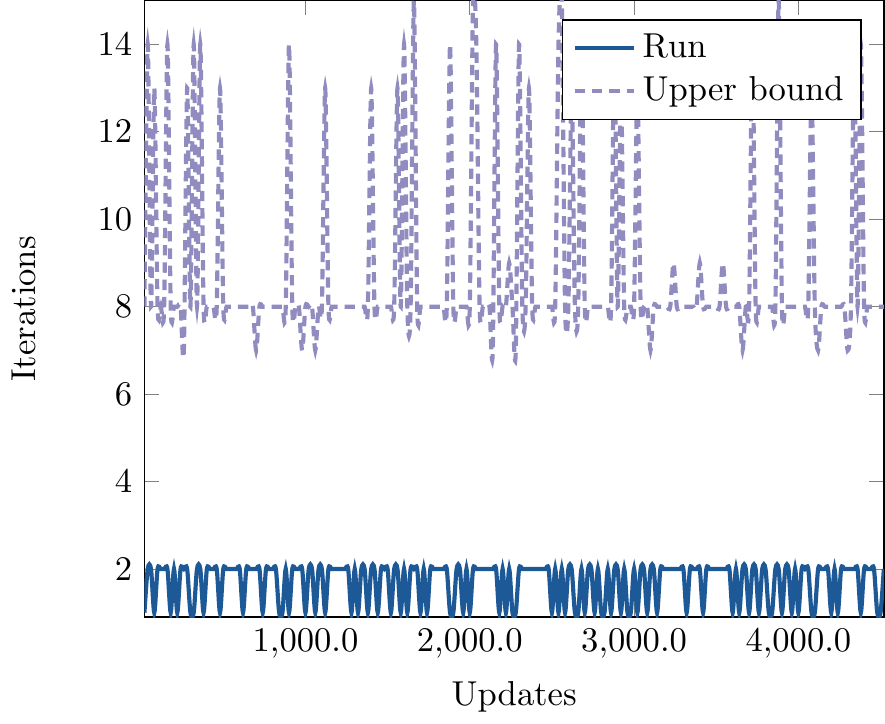}~
	    \includegraphics[width=0.32\linewidth]{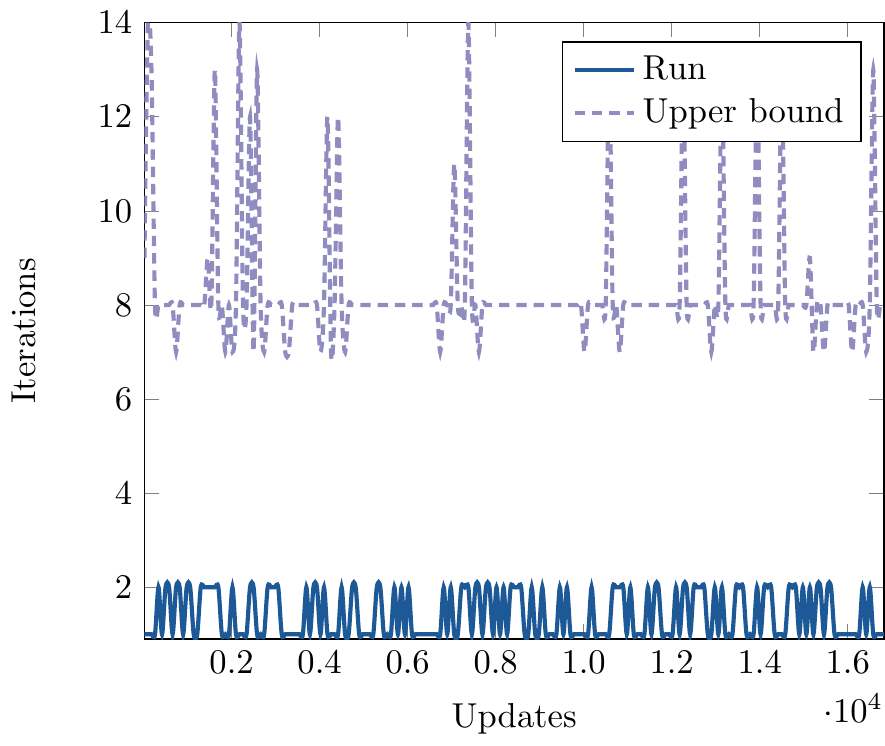}~
	    \includegraphics[width=0.32\linewidth]{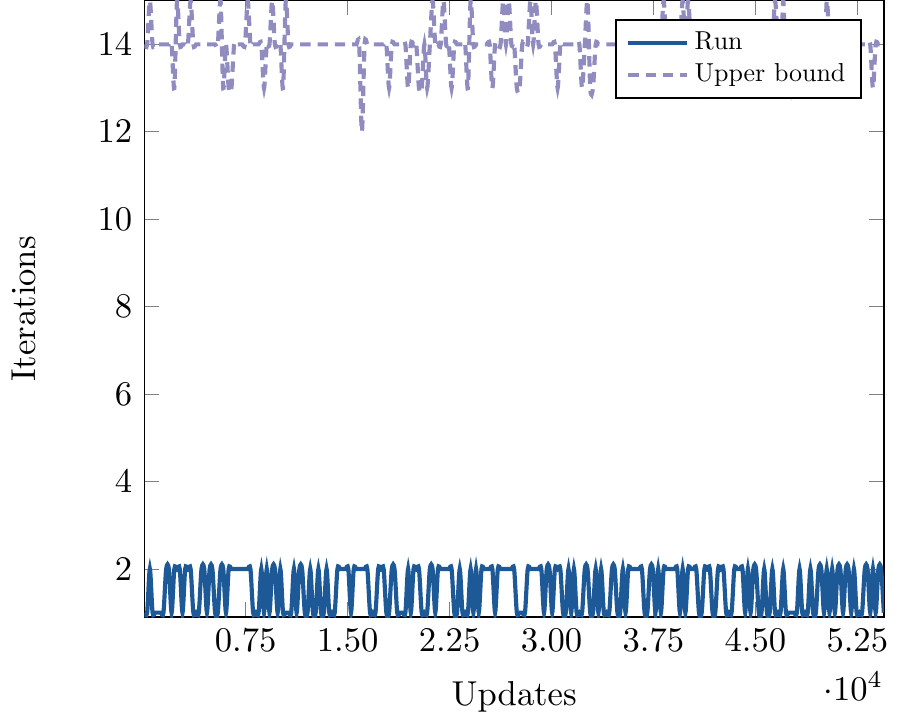}\\
    	\includegraphics[width=0.32\linewidth]{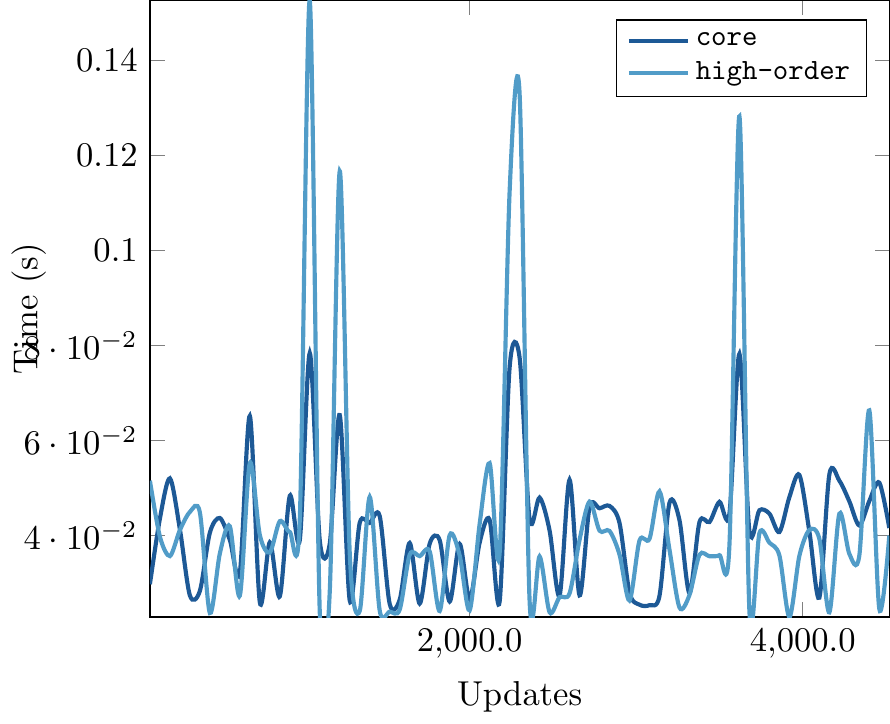}~
    	\includegraphics[width=0.32\linewidth]{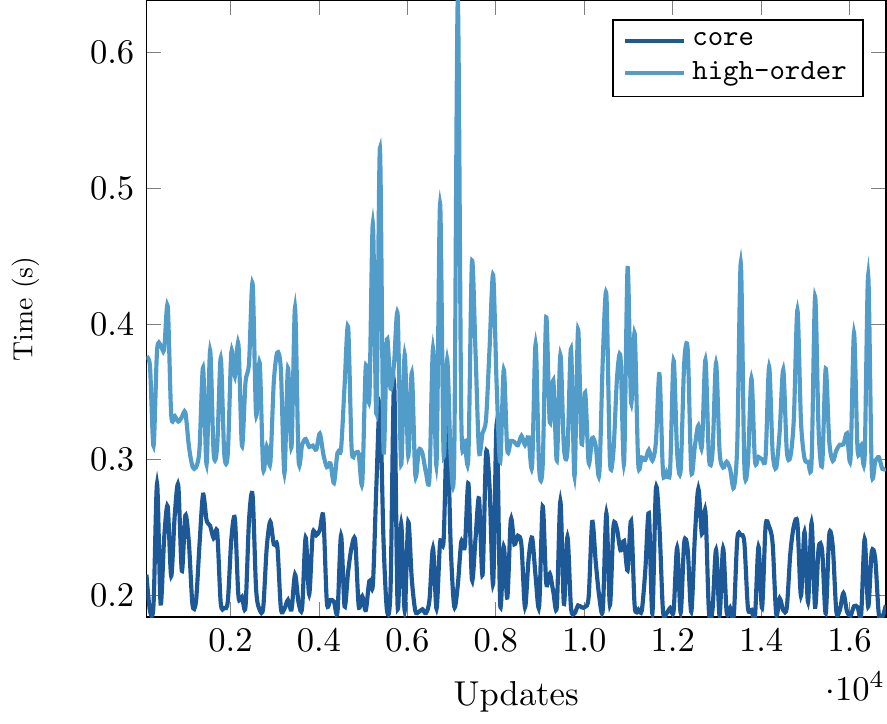}~
    	\includegraphics[width=0.32\linewidth]{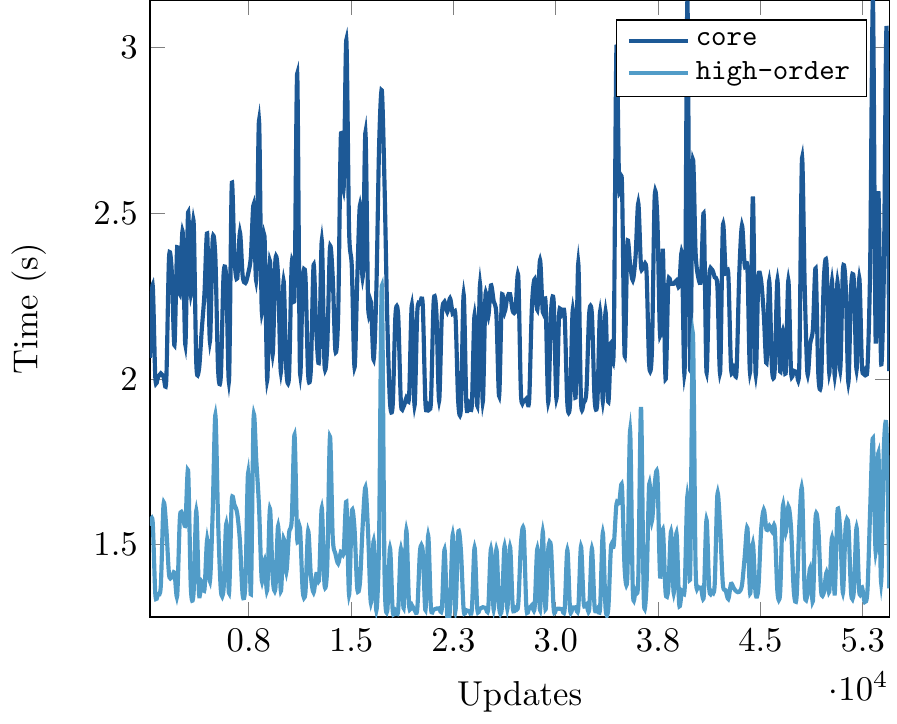}
    \caption{Performance of \cref{alg:high-level} on the \texttt{reddit-*}
    datasets, using \cref{alg:block-krylov} as the iterative method. Left to
    right: \texttt{small, medium, large} with $20$, $25$, and $50$ edge
    additions per update. Top row: number of iterations per update, which appear
    to be practically constant and dimension-independent.
    Bottom row: time elapsed computing $V$ (\texttt{core}) and higher order
    eigenvalues (\texttt{high-order}); both times are comparable to each other,
    with the total scaling linearly as $n$ increases.}
    \label{fig:reddit-datasets}
\end{figure}

Next, we evaluate our method on the \texttt{reddit-*} datasets, as depicted in
\cref{fig:reddit-datasets}. In this case, \cref{alg:block-krylov}
requires just a handful of iterations per update. Our upper bound is at most an
order of magnitude loose and independent of the underlying problem dimension.
We are thus able to handle large problems with $\cO(1)$ additional iterations
per update, in contrast to random initialization which would be expected to
scale with dimension.

We also observe that the measured wall-clock times spent updating the subspace
$V$ and the higher-order eigenvalues $\set{\hat{\lambda}_{r+1}, \dots}$ are
comparable. In fact, the majority of each step is spent on estimating
$\norm{E}_2$. Importantly, the time elapsed per iteration appears to scale
linearly as $n$ increases, which is the expected behavior given the edge density
in \cref{table:datasets} (here, edge density exactly controls the complexity of
matrix-vector multiplication).

\subsection{Principal Component Analysis} \label{sec:random-updates}
This section evaluates the performance of our adaptive method on Principal
Component Analysis (PCA)~\cite{Jolliffe02}. In this setting, we have a
data matrix $X \in \Rbb^{n \times d}$, containing $n$ points in $d$ dimensions.
For a target dimension $p$, we wish to compute $W \in \Rbb^{p \times p}$
so that $W$'s columns are the top-$p$ eigenvectors of $X^\top X$. These
columns then define a projection $X \mapsto X W \in \Rbb^{n \times p}$ which
can be helpful in exploratory data analysis, de-noising, clustering, and other
tasks. First, we show how to improve our perturbation bounds under common update models
that are applicable to incrementally updated PCA.

\subsubsection{Improved bounds under random perturbations}
\label{sec:pca-improved-bounds}
The perturbation bounds for subspaces employed in \cref{alg:high-level} can be
greatly simplified when the updates $E_{(t)}$ are random. Below, we consider
both general Gaussian random matrices, as well as sums of outer products of
Gaussian random vectors. We defer proofs of the technical results below
to \Cref{app:omitted-proofs}.

\paragraph{Random Gaussian perturbations} \label{sec:full-normal}

Suppose we are given a matrix $A := A_{(0)} \in \Rbb^{n \times n}$ (we assume
$A$ is square for simplicity, but the proof techniques extend to $A \in \Rbb^{m \times n}$).
Initially, let us assume that the perturbation matrix is i.i.d.\ Gaussian,
formalized below:
\begin{assumption} \label{asm:std-gaussian}
    For each $t$, the perturbation matrix $E := E_{(t)}$ satisfies
    \( (E)_{ij} \sim \cN\left(0, \frac{1}{n^2}\right) \), with all the
    elements independent of each other.
\end{assumption}

To derive a priori bounds for the subspace distances, we use the analogue of
Davis--Kahan for singular subspaces, due to Wedin~\cite{Wedin72}. Under
\cref{asm:std-gaussian}, we show that the distances
between singular subspaces are upper bounded and that the singular values of
$A, \tilde{A}$ are only off by a factor of $\tilde{\cO}(n^{-1})$. Both
statements hold with high probability, leading to
\cref{prop:gauss-perturb-wedin}.
\begin{proposition} \label{prop:gauss-perturb-wedin}
    Let~\cref{asm:std-gaussian} hold and let $A$, $\tilde{A} := A + E$ have SVDs
    given by $A = U \Sigma V^\top, \; \tilde{A}
    = \tilde{U} \tilde{\Sigma} \tilde{V}^\top$. Then, with probability at least
    $1 - c_1 \exp(-c_2 n) - 2n^{-\frac{\log
    n}{2}}$, \begin{align} \max\set{\norm{\sin \Theta(U_1, \tilde{U}_1)}, \norm{\sin \Theta(V_1, \tilde{V}_1)}}
    &\leq \frac{2 / \sqrt{n} + \sqrt{r} / n}{ \sigma_r - \sigma_{r+1} - (C_{A}
    + \log n) / n}, \label{eq:svdval-bound} \end{align} where $c_1, c_2$ are
    universal constants and $C_A$ is a constant only depending on $A$.
\end{proposition}

\paragraph{Low-rank random perturbations} \label{sec:lowrank-normal}

In a different setting, we observe updates to the covariance matrix
$A = X X^\top$ in which the perturbation matrix $E$ is not i.i.d.\
Gaussian but obeys a low-rank structure, as described below:
\begin{assumption} \label{asm:random-perturbation}
    For each $t$, the random perturbation $E_{(t)}$ satisfies
    \[
        E_{(t)} = \varepsilon_t z_t z_t^\top,
    \]
    where $z_t \sim \cN(0, 1/n)$ and $\prob{\varepsilon_t = 1} =
    \prob{\varepsilon_t = -1} = 1/2$ with $\varepsilon_t$ independent of $z_t$.
\end{assumption}
Suppose the above assumption holds and we want to retain information about a
$r$-dimensional subspace, $r \ll n$, over time. When we apply the a priori
bound from~\cref{thm:dav-kah}, assuming $\delta_r > 2 \norm{E}_2$, we are
reduced to bounding $\norm{E_{(t)} V_0}_2$ in the numerator, with $V_0 \in
\Rbb^{n \times r}$. Using standard tools from concentration of measure, we
can show the following:
\begin{proposition} \label{prop:lowrank-gaussian-dav-kah}
    Let \cref{asm:random-perturbation} hold (dropping the
    subscript for simplicity), and suppose that $\delta_r > 2 +
    \epsilon$. Let $\tilde{A} = A + E$ and let $V, \tilde{V}$ correspond to the
    leading $r$-dimensional subspaces of $A, \tilde{A}$ respectively.
    Then, with probability at least $1 - \expfun{-\frac{\epsilon^2 n}{32}}
    - n^{-1}$,
    \begin{align}
        \norm{\sin \Theta(V, \tilde{V})}_2
        \leq \frac{1 + \epsilon}{\delta_r} \left(\sqrt{\frac{r}{n}} +
        \sqrt{\frac{2\log n}{n}} \right).
        \label{eq:lowrank-gaussian-dav-kah}
    \end{align}
\end{proposition}
This bound can be directly applied for the estimate $d_t$ in \cref{alg:high-level}.

\subsubsection{Synthetic experiments}
We perform two experiments on synthetic data, starting from an appropriately
scaled Gaussian random matrix $X_0 \in \Rbb^{n \times d}$, using the two
variations of random updates described in~\Cref{sec:pca-improved-bounds}.
For each variation, we apply our adaptive method to keep track of a subspace of
dimension at most $\lfloor \sqrt{n} \rfloor$, which is updated
using~\cref{eq:candidate-size}. At each step, we record the true subspace
distance as well as the refined a priori bounds derived in the aforementioned
sections.

\begin{figure}[tb]
    \centerline{\includegraphics[width=\textwidth]{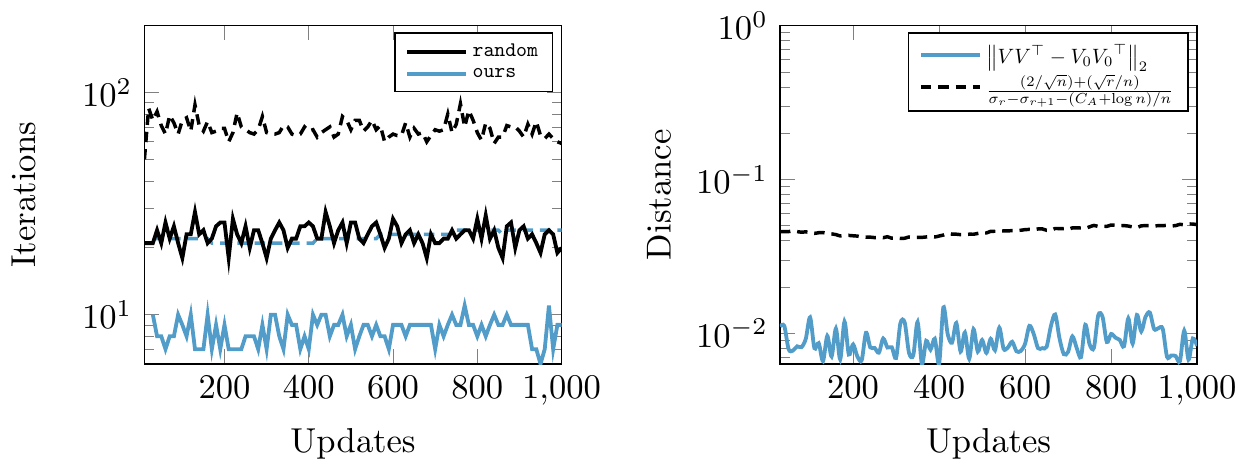}}
	\centerline{\includegraphics[width=\textwidth]{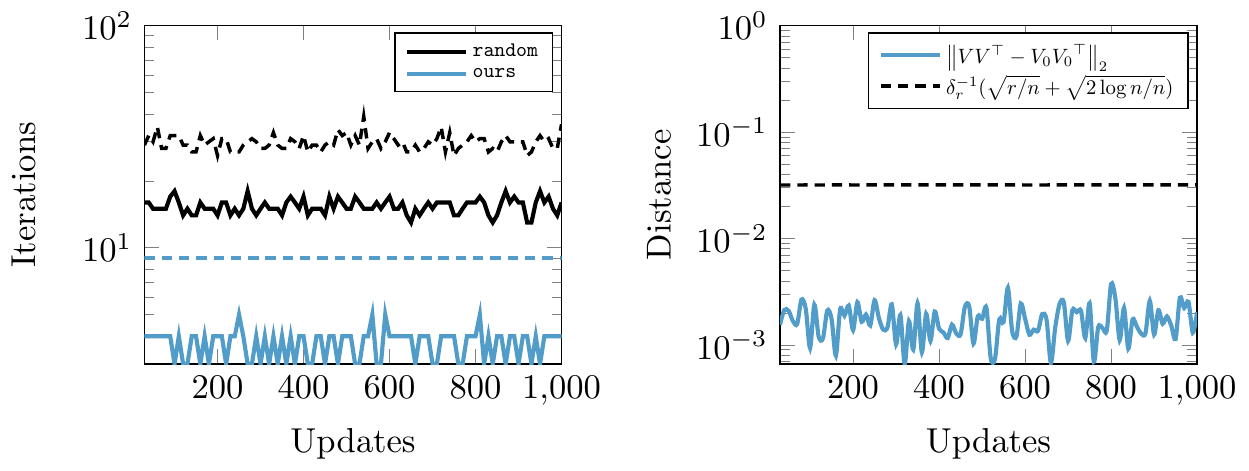}}
	\caption{PCA on synthetic data. Top: $(n,d) = (1500,1500)$
	following \cref{asm:std-gaussian}. Bottom: $(n, d) = (2000, 500)$
	following \cref{asm:random-perturbation}. Left:
    iterations of \cref{alg:block-krylov} to reach accuracy $\varepsilon =
    10^{-3}$ (solid) vs.\ upper bounds (dashed).
    Right: oracle distance vs.\ bounds
    from \cref{prop:lowrank-gaussian-dav-kah,prop:gauss-perturb-wedin}.
    Updating without warm-starting in practice is more expensive than the worst-case
    bounds of \cref{alg:high-level} even when the subspace
    distance estimates are loose.}
    \label{fig:pca-synthetic}
\end{figure}

\Cref{fig:pca-synthetic} shows the performance of \cref{alg:high-level}.
In terms of required iterations, our method clearly outperforms
reinitializing with a random matrix.
In fact, the \textit{upper bound} derived from perturbation theory falls below
the number of iterations required to reach the desired accuracy starting from a
random guess, as in~\cref{fig:prediction-comparison}. This means that
\textit{even the worst-case performance} of the warm-starting method can be
significantly better than naive seeding in simple examples.

\subsubsection{Singular Spectrum Analysis for time series data}
An application of PCA is Singular Spectrum
Analysis (SSA)~\cite{ElsTso13,GolZhi13}, which is primarily applied in
time series analysis. To apply SSA, one specifies a \textit{window length} $W$
that is expected to capture the essential behavior of the time series of length
$N \geq W$, and performs the following steps (following~\cite{GolZhi13}):
\begin{enumerate}
\item Form the trajectory matrix of the series $X$
\[
    C_X := \frac{1}{\sqrt{N}} \begin{pmatrix}
        X_{1} & X_{2} & X_{3} & \dots & X_{N - W + 1} \\
        X_{2} & X_{3} & X_{4} & \dots & X_{N - W + 2} \\
        \vdots & & & & \vdots \\
        X_{W} & X_{W+1} & X_{W+2} & \dots & X_N
    \end{pmatrix}
    \in \Rbb^{W \times (N - W + 1)}.
\]
Observe that $C_X$ is a \textit{Hankel} matrix, since its antidiagonals are
constant, and therefore admits fast matrix-vector multiplication.
\item Compute the truncated SVD $F_X$ of $C_X$ with target rank $r \ll \min(W, N - W + 1)$:
\begin{equation}
    C_X \approx F_X  = \sum_{i=1}^r u_i \sigma_i v_i^\top.
    \label{eq:cx-lr-approx}
\end{equation}
\label{step:truncated_svd}
\item Average the antidiagonals of $F_X$, from which a smoothed time series can be extracted.
\end{enumerate}

We applied SSA on the Household Power Consumption dataset, available from the
UCI Machine Learning Repository.\footnote{\url{http://archive.ics.uci.edu/ml/datasets/Individual+household+electric+power+consumption}}
The dataset contains power consumption readings for a single household spanning
4 years, spaced apart by a minute. We preprocess the dataset by computing the
moving averages of active power over $30$ minutes with a forward step of $20$
minutes, resulting in $3$ data points per hour. We set the window and series
lengths $\set{W, N} = (4032, 12096)$, which correspond to roughly 2 and 6 months
of consumption, respectively; this means that $C_X$ is the trajectory matrix
for the most recent $6$ months.

Next, we applied \cref{alg:high-level} to maintain
the decomposition of the trajectory matrix $C_X$ over time, with $r$ equal to
the value specified by the criterion in \cref{eq:candidate-size}.
Each update ``step'' moves the time series forward by $6$ measurements (or $2$
hours). The resulting dimension is determined as
$r = 4$ in every time step, satisfying
\(
    \frac{\norm{\Sigma_{1:r}}_F}{\norm{\Sigma}_F} \geq
    99\%
\)
throughout.

\begin{figure}[tb]
    \centerline{\includegraphics[width=0.9\linewidth]{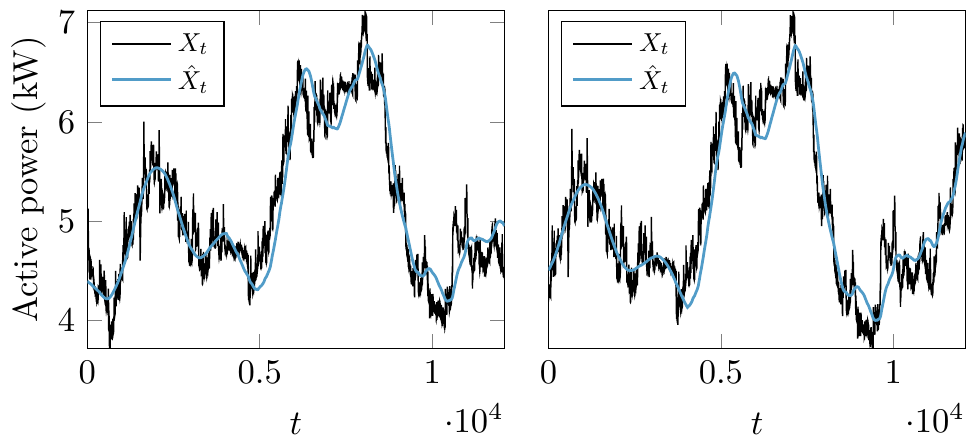}}
    \caption{Original ($X_t$) and reconstructed ($\hat{X}_t$) signals after the
    equivalent of $15$ (left) and $30$ (right) days of updates, using $r = 4$
    components and accuracy $\varepsilon = 10^{-3}$. The low-rank
    reconstruction captures the essential behavior of the time series.
    }
    \label{fig:ssa-power}
\end{figure}

\Cref{fig:ssa-power} depicts the original and reconstructed signals after $15$ and $30$ days, respectively.
The reconstructed curve captures the essential behavior of the original time
series, having ``smoothed'' out fluctuations on smaller scales. Moreover, with
the exception of the first step, all subsequent computations of the singular
subspaces take at most $3$ matrix-vector multiplications per step. On the other
hand, naive seeding can take up to $20$ matrix-vector multiplications per step.

\section{Discussion} \label{sec:discussion}

This work has detailed the theoretical and practical aspects of incrementally
updating spectral embeddings via careful use of existing information\textemdash
an oft used methodology, albeit one previously employed heuristically. Concisely,
we showed how the straightforward heuristic of warm-starting iterative eigensolvers
with previously computed subspaces can dramatically decrease
the amount of work required to update low-dimensional embeddings under small
changes to the data matrix. This is a broadly applicable scenario, encompassing
all manners of time-evolving or sequentially observed data. While we have
focused on a few specific applications, there is a wider range
of settings that fall into this regime such as spectral ranking
methods in network analysis~\cite{Perra-2008-centrality}
or latent semantic indexing in natural language processing~\cite{DDFL+90}.

A key insight, and an advantage of our proposed pipeline, is that worst-case
bounds on invariant subspace perturbations can be efficiently computed before
updates to the spectral embedding are carried out. This added information can
be used both theoretically, to bound the number of ``clean-up'' iterations that
need to be carried out, and practically, to determine if the dimension of the
subspace should change. Furthermore, the bounds simplify with additional
assumptions on the incremental updates such as sparsity or randomness.
Presenting this work in a predominantly algorithm-agnostic manner makes
the underlying conclusions applicable to any iterative method that can be
warm-started. This broadens the applicability of our analysis and strengthens
its implications.

We have focused on the \textit{dominant} eigenspace,
though \cref{alg:high-level} carries over to the complementary setting.
However, an added complexity may be the need to repeatedly solve closely related
linear systems; the nature of the incremental updates may admit ways to address
this concern, and we defer such challenges to future work. Additional challenges
pertain to removing~\cref{asm:eigval-assumption} on the spectral decay, though
we consider it rather nonrestrictive. Nevertheless, in the presence of
``challenging'' spectra, warm-starting can help reduce the cost per update and
still outperform e.g., algebraic methods, especially when matrix-vector
multiplication is cheap. Lastly, while we have considered the problem dimension
to be fixed throughout, an assumption that may be violated in practice, we
believe it is possible to pair this work with careful augmentation of the
existing subspace to address that problem.

\section*{Acknowledgements}
This research was supported by NSF Award DMS-1830274.
The authors would like to thank Yuekai Sun for his assistance in developing
the proof strategy for~\cref{prop:gauss-perturb-wedin}.

\bibliographystyle{siamplain}
\bibliography{references}

\clearpage


\appendix

\section{Auxiliary results}

\subsection{Concentration of measure}
Let us recall some results and definitions which are used in later proofs. The
following Lemma is very useful for bounding the norm of Gaussian random vectors.
It appears as part of Lemma 1 in~\cite{LauMas00}:
\begin{lemma} \label{lem:chi-squared-concentration}
    Let $(Y_1, \dots, Y_D)$ be standard Gaussian variables. Let $\alpha =
    (a_1, \dots, a_D), a_i \geq 0$ and define $Z := \sum_{i=1}^D a_i
    (Y_i^2 - 1)$.
    Then, $\forall x > 0$, we have:
    \begin{align}
        \prob{Z \geq 2 \norm{\alpha}_2 \sqrt{x} + (\max_i a_i) x}
        \leq \exp(-x)
    \end{align}
\end{lemma}

When dealing with Gaussian processes, it is common to encounter the concept of
\textit{Gaussian width}:
\begin{definition}
    Let $T \subset \Rbb^{n}$ be a bounded set and $g \sim \cN(0, I_n)$ a
    standard Gaussian random variable. The \textbf{Gaussian width} of $T$ is
    defined as
    \begin{equation}
        \mathcal{W}(T) := \expec{\sup_{x \in T} \ip{g, x}}.
        \label{eq:gauss-width}
    \end{equation}
\end{definition}

One of several bounds on the Gaussian width of a set is outlined below. It
is a straightforward consequence of~\cite[Exercise 7.6.1 \& Lemma
7.6.3]{Versh18}:
\begin{lemma} \label{lem:gauss-width-dim-bound}
    Consider $T \subset \Rbb^n$ and let $\mathrm{diam}(T) := \sup_{x \in T}
    \norm{x}_2$ and $\dim(T)$ denote its algebraic dimension. Then
    \begin{equation}
        \mathcal{W}(T) \leq \mathrm{diam}(T) \sqrt{\dim(T)}
        \label{eq:gauss-width-dim-bound}
    \end{equation}
\end{lemma}

Finally, it is known that Lipschitz functions of Gaussian variables concentrate
well around their mean. The following Theorem is standard, see
e.g.~\cite[Chapter 5]{Versh18}.
\begin{theorem} \label{thm:gaussian-lipschitz}
    Let $f: \Rbb^n \to \Rbb$ be an $L$-Lipschitz function, i.e.
    $\abs{f(x) - f(y)} \leq L \norm{x - y}$, $\forall x, y \in \Rbb^n$.
    Then, if $X = (X_1, \dots, X_n) \sim \cN(0, I_n)$, it holds that
    \begin{equation}
        \prob{f(X) - \expec{f(X)} \geq t} \leq \expfun{-\frac{t^2}{2L}}.
    \end{equation}
\end{theorem}

\subsection{Linear Algebra}

\begin{theorem}[Wedin's $\sin$ theorem] \label{thm:wedin}
    Let $A, \tilde{A} \in \Rbb^{m \times n}$ with SVDs given by
    $A = U \Sigma V^\top$, $\tilde{A} = \tilde{U} \tilde{\Sigma}
    \tilde{V}^\top$, and let $E := A - \tilde{A}$.
    Suppose there exist $\alpha, \delta > 0$ such that
    \[
        \min_{i \in [r]} \sigma_i > \alpha + \delta, \quad
        \max_{j \in [n - r]} \tilde{\sigma}_{r + j} < \alpha.
    \]
    Denote by $U_1, \tilde{U}_1, V_1, \tilde{V}_1$ the left and right singular
    subspaces corresponding to the top $r$ singular values.
    Then for any unitarily invariant norm $\norm{\cdot}$,
    \begin{equation}
        \max\set{\norm{\sin \Theta(U_1, \tilde{U}_1)},
                 \norm{\sin \Theta(V_1, \tilde{V}_1)}
        } \leq \frac{\max\set{
            \norm{E V_1}, \norm{E^\top U_1}
        }}{\delta}.
        \label{eq:wedin-sin}
    \end{equation}
\end{theorem}

\begin{theorem}[Theorem 3.6 in~\cite{ArgKnyPaiPan08}]
    \label{thm:ritz-val-distance}
    Consider a Hermitian matrix $A$ and an $A$-invariant subspace $\cX$,
    corresponding to a contiguous set of eigenvalues of $A$. Then, for any other
    subspace $\cY$ with $\dim(\cX) = \dim(\cY)$, it holds that
    \[
        \max_{i} \abs{\lambda_i(X^\top A X) - \lambda_i(Y^\top A Y)}
        \leq \rho(A) \max_i \set{ \sin^2 \theta_i (\cX, \cY) },
    \]
    where $\rho(A)$ is the spectral radius of $A$.
\end{theorem}

\begin{lemma}[Weyl's inequality] \label{lemma:weyl-ineq}
    Let $A, \widehat{A} \in \Rbb^{n_1 \times n_2}$ with $\widehat{A} = A + E$.
    Then, for all $k$, we have
    \[
        \abs{\sigma_k(A) - \sigma_k(\widehat{A})} \leq \norm{E}_2.
    \]
\end{lemma}

\begin{lemma} \label{lemma:ded-frob}
    Let $D = \mathrm{diag}(d_1, \dots, d_n)$ and $A \in \Rbb^{n \times n}$ with
    $A = A^\top$. Then it holds that
    \[
        \norm{D A D}_2 \leq \sqrt{\rank(A)} \norm{D^2 A}_2.
    \]
\end{lemma}
\begin{proof}
    Notice that we can write
    \begin{align*}
        \norm{D^2 A}_F &= \sup_{V : \norm{V}_F = 1}
        \ip{D^2 A, V} \geq \frac{\trace{A D^2 A D^2}}{\norm{D^2 A}_F}
        = \frac{\trace{D A D D A D}}{\norm{D^2 A}_F}
        = \frac{\norm{D A D}_F^2}{\norm{D^2 A}_F},
    \end{align*}
    making use of the cyclic property of the trace throughout. Rearranging
    shows that $\norm{D A D}_F \leq \norm{D^2 A}_F$. From norm equivalence
    $\norm{A}_2 \leq \norm{A}_F \leq \sqrt{\rank(A)} \norm{A}_2$, we have
    \begin{align*}
        \norm{D A D}_2 &\leq \norm{D A D}_F \leq
        \norm{D^2 A}_F \leq \sqrt{\rank(D^2 A)} \norm{D^2 A}_2,
    \end{align*}
    with $\rank(D^2 A) \leq \rank(A)$, which completes the proof.
\end{proof}

\subsection{Miscellanea}
\begin{lemma} \label{lemma:covering-number}
    Consider the set of matrices
    \[
        \cX_r = \set{X \in \Rbb^{n_1 \times n_2} \mmid
        \rank(X) \leq r, \; \norm{X}_F \leq 1}.
    \]
    Then the covering number of $\cX_r$ with respect to the metric
    $d(M_1, M_2) := \norm{M_1 - M_2}_F$ satisfies
    \[
        \cN(\cX_r, d, \varepsilon) \leq \left( \frac{9}{\varepsilon}
        \right)^{(n_1 + n_2 + 1) r}.
    \]
\end{lemma}
\begin{proof}
    Notice that we can decompose all matrices in $\cX_r$
    using the (economic) SVD as
    \(
        X = U \Sigma V^\top, \; \norm{\Sigma}_{F} \leq 1, \;
        \Sigma \in \Rbb^{r \times r}, \; U \in \Obb_{n_1, r}, \; V \in
        \Obb_{n_2, r}.
    \)
    Then we can cover the set of admissible singular values with a
    $\frac{\varepsilon}{3}$-net of vectors in $\mathbb{B}_r$, which is of
    cardinality at most
    $\left( 1 + \frac{2}{\varepsilon} \right)^{r} \leq
    \left( \frac{9}{\varepsilon} \right)^r$ (see e.g.~\cite[Chapter
    4]{Versh18}). The remainder of the proof involves covering $\Obb_{n, r}$
    and is identical to~\cite[Lemma 3.1]{CandesPlan11}.
\end{proof}

\section{Omitted proofs} \label{app:omitted-proofs}
\subsection{Proof of \cref{prop:perturb-norm-sparse-updates}}
\label{sec:sparse-updates-proof}
\newcommand{\Anew}{\tilde{A}_{\mathrm{new}}}
Given $E$, we denote $D' := \mathrm{diag}(\epsilon_1, \dots, \epsilon_n)$ where
$\epsilon_i := \sum_{j=1}^n E_{ij}$. Denoting the updated regularized adjacency
matrix by $\Anew$, it is immediate that we can write
\begin{equation}
    \Anew := (D + D')^{-1/2} (A + E) (D + D')^{-1/2}.
    \label{eq:Anew}
\end{equation}
The first step is to get an expression for $(D + D')^{-1/2}$. Since both
matrices are diagonal, we can write
\begin{align}
    (D + D')^{-1} &= \mathrm{diag}\left(\set{
        d_i\left(1 + \frac{\epsilon_i}{d_i}\right)}_{i = 1}^n
    \right)^{-1} =
    \mathrm{diag}\left(
        \set{ \frac{1}{d_i} \cdot \frac{1}{1 + \frac{\epsilon_i}{d_i}} }_{i=1}^n
    \right) \notag \\
    &= \diag\left(
        \set{ \frac{1}{d_i} \left(1 - \frac{\frac{\epsilon_i}{d_i}}{1 +
        \frac{\epsilon_i}{d_i}}\right) }_{i=1}^n
    \right) = \diag\left(\set{\frac{1}{d_i} \left(1 - \frac{\delta_i}{1 +
    \delta_i}\right)}_{i=1}^n\right) \label{eq:diagonal-deltai}, \\
    (D + D')^{-1/2} &=
    \mathrm{diag}\left(
        \set{\frac{1}{\sqrt{d_i}} \sqrt{1 - \frac{\delta_i}{1 + \delta_i}}
        }_{i=1}^n
    \right)
\end{align}
where $\delta_i := \frac{\epsilon_i}{d_i}$ in \cref{eq:diagonal-deltai}.
Next, we can write $\sqrt{1 - \frac{\delta_i}{1 + \delta_i}} = 1 - z_i$, with
$z_i \in \left[ \frac{\delta_i}{2(1 + \delta_i)}, \frac{\delta_i (2 +
\delta_i)}{2(1 + \delta_i)^2} \right]$,
based on the identity $1 - \frac{x}{2} - \frac{x^2}{2} \leq \sqrt{1 - x} \leq 1
- \frac{x}{2}$.\footnote{This identity is valid for all $x \in (-1,1)$, which
is the case for $\delta_i / (1 + \delta_i)$.}
It then follows that $(D + D')^{-1/2} = D^{-1/2} (I - Z)$,
$Z := \diag(z_1, \dots, z_n)$. Expanding this in the expression of \cref{eq:Anew}
gives us
\begin{align}
    \begin{aligned}
        \Anew &= D^{-1/2} (I - Z) (A + E) (I - Z) D^{-1/2} \\ &=
        D^{-1/2}(A - 2AZ + ZAZ)D^{-1/2}
        + (I - Z) D^{-1/2} E D^{-1/2} (I - Z) \\
        &= \tilde{A} + Z\tilde{A} Z - 2\tilde{A} Z
        + (I - Z) D^{-1/2} E D^{-1/2} (I - Z).
    \end{aligned}
    \label{eq:Anew-decomp}
\end{align}
Therefore, to bound $\norm{\Anew - \tilde{A}}_2$, we have to bound three major
components:
\[
    \norm{Z \tilde{A} Z}_2, \quad \norm{\tilde{A} Z}_2, \quad
    \text{and} \quad \norm{D^{-1/2} E D^{-1/2}}_2.
\]
For the first two terms, since $\lambda_i(\tilde{A}) \in [-1, 1]$
\cite{Chung97}, we have that:
\begin{align}
    \norm{\tilde{A} Z}_2 &\leq \norm{Z}_2 = \max_{j \in [n]} \abs{z_j}
    \leq \frac{\alpha (1 + \alpha)}{2}, \quad
    \norm{Z \tilde{A} Z}_2 \leq \left(
        \frac{\alpha(1 + \alpha)}{2} \right)^2,
    \label{eq:ZA-bound}
\end{align}
with the bound for $\abs{z_i}$ coming from the fact that
\[
    \abs{\frac{\delta_i}{2(1 + \delta_i)} + \frac{\delta_i^2}{2 (1 +
    \delta_i)^2}}
    \leq \frac{\alpha}{2} + \frac{\alpha^2}{2}, \; \text{ since }
    \abs{\delta_i} \leq \alpha.
\]
For the remaining term, we upper bound $\norm{D^{-1/2} E D^{-1/2}}_2
\leq \norm{D^{-1/2}} \norm{E D^{-1/2}}_2$, and by the Gershgorin
Circle Theorem~\cite{GVL13} (using the columns of $E D^{-1/2}$) we obtain
\begin{align}
    \begin{aligned}
    	\max_i \abs{\lambda_i(E D^{-1/2})} &\leq
    	\max_{j \in [n]} \abs{
    		\sum_{k = 1, k \neq j}^n \frac{E_{kj}}{\sqrt{d_j}}}
        \overset{(*)}{=}
    	\max_{j \in [n]}
    	\abs{\sum_{k=1}^n \frac{E_{jk}}{\sqrt{d_j}}} \leq
        \max_{j \in [n]}\frac{1}{\sqrt{d_j}} \alpha d_j \\
        \norm{D^{-1/2} E D^{-1/2}} &\leq
        \norm{D^{-1/2}}_2 \max_{j \in [n]} \alpha \sqrt{d_j} =
        \alpha \sqrt{\frac{\max_i d_i}{\min_j d_j}} = \alpha \sqrt{\kappa(D)},
    \end{aligned}
    \label{eq:ded-bound-cond}
\end{align}
where $(*)$ follows since $E_{ii} = 0, E^\top = E$ and the next inequality
follows by our assumptions. From \cref{lemma:ded-frob} (applied with $D :=
D^{-1/2}$, $A := E$) we have
\begin{align}
    \norm{D^{-1/2} E D^{-1/2}}_2 &\leq \sqrt{\rank(E)} \norm{D^{-1} E}_2
    \leq
    \sqrt{\rank(E)} \max_{i \in [n]} \abs{\sum_{j=1}^n \frac{E_{ij}}{d_i}} \leq
    \sqrt{\rank(E)} \alpha,
\end{align}
with the penultimate inequality via Gershgorin's circles. We deduce that
\begin{equation}
    \norm{D^{-1/2} E D^{-1/2}}_2 \leq
    \min\left( \sqrt{\kappa(D)}, \sqrt{\rank(E)} \right) \alpha.
    \label{eq:ded-min-bound}
\end{equation}

Combining~\cref{eq:ZA-bound,eq:ded-min-bound}
into~\cref{eq:Anew-decomp} we complete the proof (noting that $\norm{I - Z}_2^2
\leq 1 + \cO(\alpha^2)$ by submultiplicativity and the definition of $Z$):
\begin{equation}
    \norm{\Anew - \tilde{A}}_2 \leq
    \alpha \cdot \left(1 + \alpha + \sqrt{\min\set{\kappa(D), \rank(E)} }\right)
    + \left(\frac{\alpha
    (1 + \alpha)}{2}\right)^2
    \label{eq:combined-bound-normadj}
\end{equation}

\subsection{Proof of \cref{prop:gauss-perturb-wedin}}
The proof of this proposition is a combination
of \cref{lem:gauss-perturb-subspace,lemma:gauss-perturb-svdvals}. The former
controls $\norm{EV_1}_2, \norm{E^\top U_1}_2$, while the latter bounds the
difference between the corresponding singular values with high probability.
\begin{lemma} \label{lem:gauss-perturb-subspace}
    In the setting of \cref{thm:wedin}, under
    \cref{asm:std-gaussian}, it holds that
    \begin{align}
        \max \set{\norm{EV_1}_2, \norm{E^\top U_1}_2}
        \leq \frac{2}{\sqrt{n}} + \frac{\sqrt{r}}{n}
    \end{align}
    with probability at least $1 - 2\exp(-\frac{n}{2})$.
\end{lemma}
\begin{proof}
\begin{align*}
    \norm{E V_1}_2 &= \sup_{u, v \in \Sbb^{m - 1} \times \Sbb^{n - 1}}
    \ip{u, E V_1 v} =
    \sup_{u, v \in \Sbb^{m - 1} \times \mathrm{ran}(V_1)}
        \ip{u, E v}
\end{align*}
Then, appealing to the Gaussian flavor of Chevet's inequality due to
Gordon~\cite{Gordon85}, as we can write $E = \frac{1}{n} G, (G)_{ij} \sim
\cN(0, 1)$, we recover
\begin{align}
    \expec{\norm{E V_1}_2} &\leq \frac{1}{n} \left(
    \mathcal{W}(\Sbb^{n - 1}) \sup_{x \in \mathrm{ran}(V_1)} \norm{x}_2
    + \mathcal{W}(\mathrm{ran}(V_1)) \sup_{x \in \Sbb^{n-1}} \norm{x}_2
    \right)
    \leq \frac{1}{\sqrt{n}} + \frac{\sqrt{r}}{n},
    \label{eq:EV1-expec-bound}
\end{align}
where the last inequality is an appeal to \cref{lem:gauss-width-dim-bound}
for the second term after noticing $\dim(\mathrm{ran}(V_1)) = r$ combined with
the fact that
\[
    \mathcal{W}(\Sbb^{n-1}) = \expec{\sup_{v: \norm{v}_2 = 1}
    \ip{g, v}} = \expec{\norm{g}_2} \leq \sqrt{n}.
\]
An identical argument gives the same result for $\norm{E^\top U_1}_2$.

Now, using the identity $ \abs{\sup f - \sup g} \leq \abs{\sup f - g}$,
denoting $f(X) := \sup_{u, v} \frac{1}{n} \ip{u, X v}$, we recover
\begin{align*}
    \abs{f(X) - f(Y)} &\leq \frac{1}{n}
    \abs{\sup_{u, v} \ip{u, (X - Y) v}}
    \leq
    \frac{1}{n} \sup_{u,v} \cancel{\norm{u}_2} \norm{(X - Y) v}_F
    \leq \frac{1}{n} \norm{X - Y}_F,
\end{align*}
where the last two inequalities are derived using the Cauchy-Schwarz inequality.
Therefore $f$ is a $(1/n)$-Lipschitz function of Gaussian random variables.
\cref{thm:gaussian-lipschitz} gives
\begin{align}
    \prob{f(G) \geq \frac{1}{\sqrt{n}} + \frac{\sqrt{r}}{n} + t} &\leq
    \prob{f(G) \geq \expec{f} + t} \leq \exp\left( -\frac{(nt)^2}{2} \right),
    \label{eq:chevet-hp}
\end{align}
hence setting $ t = n^{-1/2}$ in \cref{eq:chevet-hp} gives us the desired high
probability bound. The proof for $\norm{E^\top U_1}_2$ is completely analogous,
and following with a union bound for the two terms we recover the desired
probability.
\end{proof}

\begin{lemma} \label{lemma:gauss-perturb-svdvals}
    In the setting of \cref{thm:wedin}, under
    \cref{asm:std-gaussian},
    it holds with probability at least $1 - 2 e^{-c_1 n} - 2n^{-\frac{\log
    n}{2}}$ that
    \begin{align}
        \abs{\sigma_{i}(\tilde{A}) - \sigma_i(A)}
        &\leq \frac{C_{A}}{n} + \frac{\log n}{n}
    \end{align}
    where $C_{A}$ is a constant depending only on $A$ and $c_1$ is
    an absolute constant.
\end{lemma}
\begin{proof}
Let us write $A = U \Sigma V^\top$ and $\hat{A} := A + E = \hat{U} \hat{\Sigma}
\hat{V}^\top$ for the singular value decompositions of the two matrices. For
a matrix $M$, write $M_{uv} := \ip{u, M v}$ for brevity and observe the
following decomposition:
\begin{align}
    \begin{aligned}
    \sigma_i(A + E) - \sigma_i(A) &= \hat{A}_{\uhat_i, \vhat_i}
    - A_{u_i, v_i}
    = A_{\uhat_i, \vhat_i} + E_{\uhat_i, \vhat_i}
    - A_{u_i, v_i} - E_{u_i, v_i} + E_{u_i, v_i} \\
    &= (E_{\uhat_i, \vhat_i} - E_{u_i, v_i}) + E_{u_i, v_i}
    + \left( A_{\uhat_i - u_i, v_i} + A_{u_i, \vhat_i - v_i} + A_{\uhat_i - u_i,
    \vhat_i - v_i} \right) \\
    &= (E_{\uhat_i, \vhat_i} - E_{u_i, v_i}) + E_{u_i, v_i}
    - \sigma_i(A) \left[ (1 - \ip{\uhat_i, u_i}^2) + (1 - \ip{\vhat_i, v_i}^2)
    \right]  \\
    &\quad + A_{\uhat_i - u_i, \vhat_i - v_i},
    \end{aligned}
    \label{eq:svdvals-decomp}
\end{align}
with the last equality above following since $A v_i = \sum_{j=1}^n \sigma_j(A)
u_j v_j^\top v_i = \sigma_i(A) u_i$.
The decomposition above hints towards the objects that we need to control to
prove that $\sigma_i(\hat{A}) - \sigma_i(A)$ is small in magnitude. To be
precise, let us define the set
\begin{equation}
    \mathbf{S}_{\delta} := \set{(u, v) \in \Sbb^{n-1} \times \Sbb^{n-1}
    \mmid \norm{uv^\top - u_i v_i^\top}_F^2 \leq
    \delta^2, \; \ip{u, u_i} \geq 0, \; \ip{v, v_i} \geq 0}
    \label{eq:sdelta}
\end{equation}
and the events below:
\begin{align}
    \begin{aligned}
    \cE_1 &:= \set{ \norm{E}_2^2 \leq \frac{C_1(x)}{n}},
    \quad
    \cE_2 := \set{\abs{E_{u_i, v_i}} \leq \frac{C_2(x)}{n}}, \\
    \cE_3 &:= \set{\sup_{u, v \in \mathbf{S}_{\delta}} \abs{E_{u, v}
    - E_{u_i, v_i}} \leq \frac{\delta}{\sqrt{n}} + \frac{\delta
    C_3(x)}{n}},
    \end{aligned}
    \label{eq:gauss-svd-highprob-events}
\end{align}
where $C_1, C_2, C_3$ only depend on $x$. The set $\mathbf{S}_{\delta}$ contains
the unit vectors $(u, v)$ such that $E_{u, v}$ is ``close'' to $E_{u_i, v_i}$
in terms of $L^2$ distance, with the condition $\ip{u, u_i}, \ip{v, v_i} \geq 0$
verifiable from the variational characterization of singular values. The goal
is to control the probability of each of $\cE_i$ being false.

As a first step, notice that~\cref{thm:wedin} implies
\begin{equation}
    \norm{u_i u_i^\top - \hat{u}_i \hat{u}_i^\top}_F^2 +
    \norm{v_i v_i^\top - \hat{v}_i \hat{v}_i^\top}_F^2
    \leq \frac{2}{\delta_{\mathrm{Wedin}}^2} \norm{E}_2^2
    \leq \frac{2}{\left(\alpha - \frac{1}{\sqrt{n}}\right)^2} \frac{C}{n}
    = \frac{c_1}{n},
    \label{eq:svddist-wedin-bound}
\end{equation}
assuming that $\alpha > \frac{1}{\sqrt{n}}$, with the last inequality above
holding with high probability since $\norm{E}_2 \asymp \frac{1}{\sqrt{n}}$ from
elementary arguments in random matrix theory. Notice that the above in
particular is equivalent to
\begin{align}
    2 (1 - \ip{u_i, \hat{u}_i}) + 2
    (1 - \ip{v_i, \hat{v}_i}) \leq \frac{c_1}{n}, \quad
    \max \set{1 - \ip{u_i, \hat{u}_i}, 1 - \ip{v_i, \hat{v}_i}}
    \leq \frac{c_1}{n}
\end{align}
since $1 - \ip{u_i, \hat{u_i}} \leq 1 - \ip{u_i, \hat{u}_i}^2$ due to all
vectors being unitary and $\ip{u_i, \hat{u}_i} \geq 0$; the same line of
reasoning applies to $\ip{v_i, \hat{v}_i}$.
Equivalently, the above gives us
\begin{equation}
    \norm{u_i - \hat{u}_i}_2^2 + \norm{v_i - \hat{v}_i}_2^2 \leq \frac{c_1}{n},
    \quad \max\set{\norm{u_i - \hat{u}_i}_2, \norm{v_i - \hat{v}_i}_2}
    \leq \frac{c_2}{\sqrt{n}}
    \label{eq:delta_range}
\end{equation}
The magnitude of the RHS in \cref{eq:delta_range} will help us determine the
correct range for $\delta$ in $\cE_3$.

For the remainder, consider the Gaussian process $\set{n E_{u, v} - n
E_{u_i, v_i}}_{(u, v) \in \Sbb^{n-1} \times \Sbb^{n-1}}$. Following arguments
from~\cite[Chapter 7]{Versh18}, since $nE$ is equal in distribution to a
standard Gaussian matrix, we obtain using the standard $L^2$ distance:
\begin{equation}
    d((u, v), (u', v'))^2 :=
    \norm{n E_{u,v} - n E_{u', v'}}^2_{L^2}
    = \norm{u v^\top - u' {v'}^\top}_F^2
    \leq \norm{u - u'}_2^2 + \norm{v - v'}_2^2,
    \label{eq:gaussian-l2-dist}
\end{equation}
where the inequality follows from~\cite[Exercise 7.3.2]{Versh18}.
We can now define the following quantities, which control most of the terms
above:
\begin{equation}
    \xi_{\delta} := \sup_{u, v \in \mathbf{S}_{\delta}} E_{u, v} - E_{u_i,
    v_i}, \quad
    \tilde{\xi}_{\delta} := \sup_{u, v \in \mathbf{S}_{\delta}}
    \abs{E_{u, v} - E_{u_i, v_i}}.
    \label{eq:gaussian-xi-delta}
\end{equation}
Since $0 \in \mathbf{S}_{\delta}$ by setting $(u, v) = (u_i, v_i)$, it follows
that
\[
    \sup_{u, v \in \mathbf{S}_{\delta}} E_{u, v} - E_{u_i, v_i}, \quad
    \sup_{u, v \in \mathbf{S}_{\delta}} -(E_{u, v} - E_{u_i, v_i})
    \geq 0
\]
and also that $\norm{\tilde{\xi}_{\delta}}_{\psi_2} \asymp
\norm{\xi_{\delta}}_{\psi_2}$, with $\norm{\cdot}_{\psi_2}$ denoting the
subgaussian norm. To arrive at the result, we need the following Lemma:
\begin{lemma} \label{lemma:gaussian-variance-proxy}
    For $\tilde{\xi}_{\delta}$ defined as in \cref{eq:gaussian-xi-delta}, it
    holds that
    \[
        \norm{\tilde{\xi}_{\delta}}_{\psi_2}
        \lesssim \left( \frac{\delta}{n} \right)^2, \;
        \quad \expec{\tilde{\xi}_{\delta}} \lesssim \frac{\delta}{\sqrt{n}}.
    \]
\end{lemma}
\begin{proof}
    The proof of the first property is immediate
    since~\cref{eq:gaussian-l2-dist} and the fact that we are working within
    $\mathbf{S}_{\delta}$ give us
    \begin{align*}
        \sup_{u, v \in \mathbf{S}_{\delta}} \norm{E_{u, v} - E_{u_i,
        v_i}}^2_{L^2}
        = \sup_{u, v \in \mathbf{S}_{\delta}}
        \frac{1}{n^2}
            \norm{uv^\top - u_i v_i^\top}_F^2
        \leq \frac{\delta^2}{n^2},
    \end{align*}
    hence the result for the subgaussian norm follows from the high-probability
    version of Dudley's theorem~\cite[Theorem 8.1.6]{Versh18}.

    For the latter property, let us write $Z_{u, v} = n E_{u, v}$, so that
    $Z := n E$ is a standard Gaussian matrix. Then, using Dudley's inequality,
    it follows that
    \begin{align}
        \expec{\xi_{\delta}} &= \frac{1}{n} \expec{\sup_{u, v \in
        \mathbf{S}_{\delta}} Z_{u, v} - Z_{u_i, v_i}}
        \leq \frac{c_1}{n} \int_0^{\delta} \sqrt{\log \cN(\mathbf{S}_{\delta},
        d, \varepsilon)} \dd{\varepsilon},
        \label{eq:proof-dudley-integral}
    \end{align}
    with $\cN(\cS, d, \cdot)$ denoting the \textit{covering number} of $\cS$
    with respect to the metric $d$, and the upper limit
    $\delta$ following since $\cN(\mathbf{S}_{\delta}, d, \varepsilon)
    = 1, \; \forall \varepsilon > \delta$ as the distance between two ``points''
    indexed by $\mathbf{S}_{\delta}$ is never more than $\delta$.

    To estimate $\cN(\mathbf{S}_{\delta}, d, \varepsilon)$, it suffices to
    notice that
    \[
        d((u, v), (u_i, v_i)) = \norm{u v^\top - u_i v_i^\top}_F
        \in \set{X \in \Rbb^{n \times n} \mmid \rank(X) \leq 2, \;
        \norm{X}_F \leq \delta},
    \]
    which is just a scaled version of $\cX_2$ appearing
    in~\cref{lemma:covering-number}. Therefore it must satisfy
    \[
        \cN(\mathbf{S}_{\delta}, d, \varepsilon) \leq
        \cN(\delta \cX_2, d, \varepsilon) \lesssim \left(\frac{9
        \delta}{\varepsilon}\right) ^{4n + 2}.
    \]
    Substituting the above into the integral
    of \cref{eq:proof-dudley-integral}, we recover that
    \[
        \expec{\xi_{\delta}} \leq \frac{c_1}{n}
        \int_0^{\delta} \sqrt{(4n+2) \log \frac{9\delta}{\varepsilon}}
        \dd{\varepsilon} \asymp \frac{\delta}{\sqrt{n}}.
    \]
\end{proof}
Given~\cref{lemma:gaussian-variance-proxy}, we can deduce that
\begin{align*}
    \prob{\tilde{\xi}_{\delta} \geq c_1 \frac{\delta}{\sqrt{n}} + t} \leq
    \expfun{-c_2\left( \frac{n}{\delta} \right)^2 t^2} \Rightarrow
    \prob{\tilde{\xi}_{\delta} \geq 2 \frac{\delta}{\sqrt{n}}}
    \leq \expfun{-c_2 n}.
\end{align*}
Moreover, we know that $\ip{u_i, E v_i} \sim \frac{1}{n} g, \; g \sim \cN(0, 1)$
and therefore it follows from Gaussian concentration that
\[
    \prob{\abs{E_{u_i, v_i}} \geq t}
    \leq 2 \expfun{-\frac{n^2 t^2}{2}} \Rightarrow
    \prob{\abs{E_{u_i, v_i}} \geq \frac{\log n}{n}} \leq 2 n^{-\frac{\log
    n}{2}}.
\]
Finally, an appeal to~\cite[Corollary 7.3.3]{Versh18} gives us that
\[
    \prob{\norm{E}_2 \geq \frac{2 + c_3}{\sqrt{n}}} \leq 2 \expfun{-c_4 n},
\]
so we can deduce (via a union bound) that $\prob{\setI_{i=1}^3 \cE_i}
\geq 1 - 2 \exp(-c n) - 2 n^{-\frac{\log n}{2}}$, by setting $c := \max(c_2,
c_4)$ using the notation above. Then, returning to~\cref{eq:svdvals-decomp},
we deduce that
\begin{align*}
    \abs{\sigma_i(A + E) - \sigma_i(A)} &\leq
    2 \frac{\delta}{\sqrt{n}} + \frac{\log n}{n} +
    \sigma_i(A) \frac{C}{n} + \norm{A}_2 \frac{C'}{n}
\end{align*}
Looking back at \cref{eq:delta_range}, it is immediate that $\delta \lesssim
\frac{1}{\sqrt{n}}$, which completes the desired claim.
\end{proof}

Putting \Cref{lem:gauss-perturb-subspace,lemma:gauss-perturb-svdvals} together
and appealing to Weyl's inequality for the singular value gap, we deduce that
under \cref{asm:std-gaussian} the desired inequality must hold
with probability at least $1 - c_1 \exp(-c_2 n) - 2n^{-\frac{\log n}{2}}$.

\subsection{Proof of \cref{prop:lowrank-gaussian-dav-kah}}
The proof consists of two components. In order for \cref{thm:dav-kah}
to be applicable, we need to ensure that $\norm{E}_2 < \frac{\delta_r}{2}$, so
given the lower bound on $\delta_r$ it would suffice to prove that
\(
    \norm{E}_2 = \norm{\varepsilon_t z_t z_t^\top}_2
    = \norm{z_t}^2 \leq 1 + \frac{\epsilon}{2},
\)
which is a straightforward application of $\chi^2$-concentration. The other
component consists of showing that $\norm{E V_0}_2$ is small with high
probability, which we tackle by bounding its expectation and using Gaussian
concentration to obtain a high probability bound.

Let us take care of the first component; we can write
$\norm{z_t}_2^2 = \sum_{i=1}^n \frac{1}{n} g_i^2, \; g_i \sim \cN(0, 1)$.
Then \cref{lem:chi-squared-concentration} with $\alpha = \frac{1}{n}
\bm{1}$ tells us that
\begin{align*}
    \prob{\norm{E}_2 \geq 1 + 2 \sqrt{\frac{x}{n}} + \frac{x}{n}}
    \leq \exp(-x)
\end{align*}
Setting $x = \frac{\epsilon^2 n}{32}$ when $\epsilon$ is small gives us that
$\norm{E}_2 \leq 1 + \frac{\epsilon}{2 \sqrt{2}} + \frac{\epsilon^2}{32}
< 1 + \frac{\epsilon}{2}$ with probability at least $1 -
\expfun{-\frac{\epsilon^2 n}{32}}$.

In order to address the second component, we prove the following Lemma:
\begin{lemma} \label{lemma:rand-perturbation-davis-kahan-numerator}
    Suppose~\cref{asm:random-perturbation} holds. Then, for any $r$-dimensional
    subspace $V_0$, there exist constants $c_1, c_2 > 0$ such that
    \begin{align*}
        \prob{\norm{E^{(t)} V_0}_2 \leq c_1 \sqrt{\frac{r \log n}{n}}} \geq 1 -
        \expfun{-c_2 n} - n^{-r}.
    \end{align*}
\end{lemma}
\begin{proof}
    Denote $E \equiv E^{(t)}$ for brevity.
    Using the variational characterization of singular values, we know that
    \begin{align*}
        \norm{E V_0}_2 &= \sup_{u \in \Sbb^{n - 1}, v \in \Sbb^{r - 1}}
            \ip{u, E V_0 v}
            = \sup_{u \in \Sbb^{n - 1}, v \in \Sbb^{r - 1}}
            \varepsilon_t \ip{u, z_t z_t^\top V_0 v}
            \leq \norm{z_t}_2 \sup_{v \in V_0 \Sbb^{r - 1}}
            \ip{z_t, v}.
    \end{align*}
    Taking expectations and applying the Cauchy-Schwarz inequality:
    \begin{align*}
        \expec{\norm{E V_0}_2} &\leq
        \expec{\norm{z_t}^2}^{1/2}
        \expec{\sup_{v \in V_0 \Sbb^{r - 1}} \ip{z_t,
        v}^2}^{1/2}
        \leq \left( \frac{1}{n}
        \cW(\Sbb^{n - 1} \cap \mathrm{ran}(V_0))^2 \right)^{1/2},
    \end{align*}
    where the last inequality follows from the fact that $\expec{\norm{z_t}^2}
    = \frac{1}{n}$, the definition of gaussian width, and the fact that $V_0$
    is an orthogonal matrix, hence $V_0$ acting on the unit sphere maps back to
    the unit sphere. However, $\mathrm{diam}(\Sbb^{n-1} \cap \mathrm{ran}(V_0))
    \leq 1$, and $\dim(T) = \dim(\mathrm{ran}(V_0)) = r$. An appeal to
    \cref{lem:gauss-width-dim-bound} recovers
    \(
        \expec{\norm{E V_0}_2} \leq \sqrt{\frac{r}{n}}.
    \)

    We can now proceed to show that $\norm{E V_0}_2$ is small with high
    probability. Let us denote the events
    \[
        \cE_1 := \set{\norm{z_t}_2 \leq c_1}, \quad
        \cE_2 := \set{\sup_{v \in V_0 \Sbb^{r - 1}} \ip{z_t, v} \leq
        \sqrt{\frac{r}{n}} +
        \sqrt{\frac{2 r \log n}{n}}}.
    \]
    For the first event, we know from \cref{lem:chi-squared-concentration}
    that $\prob{\cE_1^c} \leq \exp(-c_2 n)$, where $c_2$ depends only on $c_1$.
    Additionally, writing $z_t = \frac{1}{\sqrt{n}} g, \; g \sim \cN(0, I_n)$,
    we see that
    \begin{align*}
        \abs{\sup_{v \in V_0 \Sbb^{r - 1}} \ip{z_t, v}
        - \sup_{v \in V_0 \Sbb^{r - 1}} \ip{z'_t, v}}
        &\leq \abs{\sup_{v \in V_0 \Sbb^{r - 1}}
        \ip{z_t - z'_t, v}} \leq
        \frac{1}{\sqrt{n}} \norm{g - g'},
    \end{align*}
    which means that $f(X) := \frac{1}{\sqrt{n}} \sup_{v \in V_0 \Sbb^{r - 1}}
    \ip{X, v}$ is a $(1 / \sqrt{n})$-Lipschitz function of Gaussian variables.
    \Cref{thm:gaussian-lipschitz} then implies
    \begin{align*}
        \prob{\sup_{v \in V_0 \Sbb^{r - 1}} \ip{z_t, v}
        \geq \expec{\sup_{v \in V_0 \Sbb^{r - 1}} \ip{z_t, v}} + t}
        &\leq \exp\left( -\frac{t^2}{2 L^2} \right) \\
    \Rightarrow
    \prob{\cE_2^c}
    = \prob{\sup_{v \in V_0 \Sbb^{r - 1}} \ip{z_t, v}
        \geq \sqrt{\frac{r}{n}} + \sqrt{\frac{2 \log n}{n}}}
        &\leq \frac{1}{n}.
    \end{align*}
    Taking a union bound over $\cE_1^c, \cE_2^c$ we recover
    \[
        \prob{\norm{E V_0}_2 \leq c_1 \sqrt{\frac{r}{n}}
        + c_1 \sqrt{2} \sqrt{\frac{\log n}{n}}
        } \geq 1 - \expfun{-c_2 n} - n^{-1}.
    \]
\end{proof}
Finally, to conclude the proof, we set $c_1 = 1 + \frac{\epsilon}{2 \sqrt{2}}
+ \frac{\epsilon^2}{32}$ in the definition of $\cE_1$ above, which gives $c_2 =
\frac{\epsilon^2}{32}$. \qed

\end{document}